\documentclass[11pt]{amsart}
\usepackage[a4paper,top=0.8in,bottom=0.8in,left=0.8in,right=0.8in]{geometry}
\usepackage{libertinus}
\usepackage{setspace}
\usepackage{mathtools}
\usepackage{eurosym,color,scalefnt}
\usepackage{amsmath,amsfonts,float,array}
\usepackage{graphicx,amssymb}
\usepackage{array}
\usepackage{longtable}
\usepackage{booktabs}
\usepackage{geometry}
\usepackage{placeins,graphicx,float,xcolor,tikz}
\newcolumntype{S}{>{\raggedright\arraybackslash}p{0.80\textwidth}}
\newcolumntype{J}{>{\raggedright\arraybackslash}p{0.30\textwidth}}

\newenvironment{prooftable}
{\noindent\renewcommand{\arraystretch}{1.3}%
	\begin{longtable}{r S J}%
		\hline
		\textbf{Step} & \textbf{Statement} & \textbf{Justification} \\
		\hline
		\endfirsthead
		\hline
		\textbf{Step} & \textbf{Statement} & \textbf{Justification} \\
		\hline
		\endhead
		\hline
		\multicolumn{3}{r}{\textit{Continued on next page}} \\
		\endfoot
		\hline
		\endlastfoot}
	{\end{longtable}}
\usepackage[utf8]{inputenc}
\usepackage[colorlinks=true,
citecolor=red,
urlcolor=blue]{hyperref}
\usepackage{ragged2e}
\usepackage{xcolor}
\usepackage{titlesec,float}
\usepackage{tikz}
\usepackage{pgfplots}

\usepackage{fancyhdr}
\usepackage{amsmath, amssymb, amsthm}
\usepackage{stmaryrd}
\usepackage{setspace}
\usepackage{listings}
\usepackage{xcolor}
\usepackage{logicproof}
\usepackage{array}
\usepackage{graphicx}
\usepackage{float}
\usepackage{adjustbox}
\usepackage{titlesec}
\usepackage{titletoc}
\usepackage{fancyhdr}
\usepackage{lastpage}
\usepackage[colorlinks=true,
citecolor=red,
urlcolor=blue]{hyperref}

\usepackage{amsmath, amssymb}
\usepackage{enumitem}
\usepackage{mathrsfs}  
\usepackage{hyperref}
\usepackage{xcolor}

\definecolor{darkblue}{RGB}{0,0,128}
\definecolor{darkgreen}{RGB}{0,100,0}
\definecolor{darkred}{RGB}{128,0,0}
\definecolor{purple}{RGB}{128,0,128}
\definecolor{orange}{RGB}{255,140,0}

\newcommand{\lo}{\lowercase} 

\newtheorem{theorem}{Theorem}[section]
\newtheorem{corollary}{Corollary}[section]
\newtheorem{lemma}{Lemma}[section]
\newtheorem{proposition}{Proposition}[section]

\newtheorem{algorithm}{Algorithm}[section]
\theoremstyle{definition}
\newtheorem{definition}{Definition}[section]

\newtheorem{example}{Example}[section]
\newtheorem{notation}{Notation}[section]
\theoremstyle{remark}
\newtheorem{remark}{Remark}[section]

\numberwithin{equation}{section}

\definecolor{darkblue}{RGB}{0,0,128}
\definecolor{darkgreen}{RGB}{0,100,0}
\definecolor{darkred}{RGB}{128,0,0}
\definecolor{purple}{RGB}{128,0,128}
\definecolor{orange}{RGB}{255,140,0}
\makeatletter
\@namedef{subjclassname@2020}{
	\textup{2025} Mathematics Subject Classification
}
\makeatother

\titleformat{\section}
{\normalfont\Large\bfseries\color{blue}}
{\thesection}{1em}{}
\titleformat{\subsection}
{\normalfont\large\bfseries\color{teal}}
{\thesubsection}{1em}{}
\titleformat{\subsubsection}
{\normalfont\normalsize\bfseries\color{violet}}
{\thesubsubsection}{1em}{}
\titlespacing*{\section}{0pt}{12pt}{6pt}
\titlespacing*{\subsection}{10pt}{10pt}{4pt}
\titlespacing*{\subsubsection}{0pt}{8pt}{2pt}
\makeatletter
\renewcommand{\l@section}{\@tocline{1}{1em}{1em}{2em}{}}
\renewcommand{\l@subsection}{\@tocline{1}{1em}{3em}{2em}{}}
\makeatletter
\title{M\lo{odular} r\lo{esolutions by} p\lo{olyseries}}
\author{  B\lo{rahimi}, M\lo{ahdi}. T\lo{ahar}. \\
	{{D\lowercase{epartment of} M\lowercase{athematics}}} \\
	{U\lowercase{niversity of} M\lowercase{ohamed} B\lowercase{oudiaf}, M\lowercase{sila}, A\lowercase{lgeria}} \\
	{ \texttt{\href{mailto:mahdi.brahimi@univ-msila.dz}{%
				\lowercase{	mahditahar.brahimi@univ-msila.dz}}}} }
\makeatother
\usetikzlibrary{positioning, shapes.geometric, arrows.meta}
\tikzset{
	algebra/.style={
		rectangle,
		rounded corners,
		draw=black,
		fill=white,
		text width=3cm,
		align=center,
		font=\small
	},
	wide/.style={
		algebra,
		text width=4cm
	},
	arr/.style={-Stealth, thick}
}

\begin{document}
\setstretch{1}\scriptsize	
\begin{abstract}
	We study the modular resolution method using new tools called polynumbers
	and polyseries, introduced by Prof. Wildberger N.J. We try to prove an
	equivalence theorem of the existence and the uniqueness of the solutions of
	the modular quadratic equations, using the recurrence formula between the
	Catalan sequence terms and introducing the following notions: Wildberger's
	polynumber sequences (polynomials), binomial Chu-Vandermonde identity and
	truncated polyseries.
\end{abstract}
\maketitle
\thispagestyle{empty}
\tableofcontents
\newpage
\pagestyle{plain}
 \pagenumbering{arabic} \vspace{0.3in}

\section{Discretization}

\pagenumbering{arabic} Modular resolutions
help us to understand the structure of modules over rings with positive
characteristic. In this paper, we introduce a new technique based on
Wildberger's polyseries expansions to construct and analyze such resolutions.

\begin{remark}
	\text{ \ \ }\newline
	Most of the present results are indebted to the wise efforts of Prof.
	Wildberger, N.J. \newline
	For more details, \textbf{See} \cite%
	{njwildbergerMathFoundations,Wildberger2015Data, WildEggmathematicscourses,
		Wildberger2021Box, Wildberger2021Poly, njwildbergerAlgCal}
\end{remark}

\subsection{Polynumber Sequences and polyseries}

\begin{definition}
	A polynumber $u=\left[ u\left( i\right) \right\rangle _{0}^{m}\equiv
	\sum_{i=0}^{m}a_{i}e_{i}$ , \ $m\in 
	\mathbb{N}
	,$ is a vexel which supports congruent addition and multiplication, the
	multiset extension of the following congruences identities: 
	\begin{equation*}
	\alpha ^{n}\,.\,\alpha ^{m}\equiv e_{n}\,.\,e_{m}\equiv \left[ \left[ n%
	\right] \right] \,.\,\left[ \left[ m\right] \right] \equiv e_{n+m}\equiv
	\alpha ^{n+m}\equiv \left[ \left[ n+m\right] \right] \qquad n,\text{ }m\in 
	\mathbb{N}
	\quad 
	\end{equation*}%
	where: $\ \ e_{i}=\left[ \delta _{i,j}\right\rangle _{j\in 
		\mathbb{N}
	},$ is defined by the Kronecker delta symbol:%
	\begin{equation*}
	\ \delta _{i,j}=\left\{ 
	\begin{array}{ccc}
	1 &  & i=j \\ 
	&  &  \\ 
	0 &  & i\neq j%
	\end{array}%
	\right. 
	\end{equation*}%
	
	For two polynumbers $\ u\equiv \sum_{i=0}^{m}a_{i}e_{i}$ and $v\equiv
	\sum_{j=0}^{n}b_{j}e_{j}$ we have:
	\begin{equation*}
	u+v\equiv \sum_{k=0}^{\max (m,n)}(a_{k}+b_{k})e_{k}
	\end{equation*}
	
	and \ 
	\begin{equation*}
	u\cdot v\equiv \sum_{k=0}^{m+n}\left( \sum_{i+j=k}a_{i}b_{j}\right) e_{k}
	\end{equation*}
\end{definition}

\begin{example}
	\text{ \ \ }
	
	\begin{itemize}
		\item \text{ \ \ }If $\ u\equiv \left[ \left[ 0\right] \text{ }\left[ 1%
		\right] \text{ }\left[ 3\right] \text{ }\left[ 1\right] \right] $, then 
		\begin{equation*}
		u\equiv e_{0}+2e_{1}+e_{3}\equiv 
		\begin{tabular}{|l|}
		\hline
		$1$ \\ 
		$2$ \\ 
		$0$ \\ 
		$1$%
		\end{tabular}%
		,
		\end{equation*}%
		and if $\ v\equiv \left[ \left[ 2\right] \text{ }\left[ 2\right] \text{ }%
		\left[ 1\right] \text{ }\right] ,$ then 
		\begin{equation*}
		v\equiv e_{1}+2e_{2}\equiv 
		\begin{tabular}{|l|}
		\hline
		$0$ \\ 
		$1$ \\ 
		$2$%
		\end{tabular}%
		,
		\end{equation*}%
		such that:%
		\begin{eqnarray*}
			u\,.\,v &\equiv &\left[ \left[ 0\right] \text{ }\left[ 1\right] \text{ }%
			\left[ 3\right] \text{ }\left[ 1\right] \right] \,.\,\left[ \left[ 2\right] 
			\text{ }\left[ 2\right] \text{ }\left[ 1\right] \text{ }\right]  \\
			&\equiv &\left[ \left[ 2\right] \text{ }\left[ 2\right] \text{ }\left[ 1%
			\right] \text{ }\left[ 3\right] \text{ }\left[ 3\right] \text{ }\left[ 2%
			\right] \text{ }\left[ 5\right] \text{ }\left[ 5\right] \text{ }\left[ 4%
			\right] \text{ }\left[ 3\right] \text{ }\left[ 3\right] \text{ }\left[ 2%
			\right] \right]  \\
			&\equiv &\left( e_{0}+2e_{1}+e_{3}\right) \left( e_{1}+2e_{2}\right)
			=e_{1}+4e_{2}+4e_{3}+e_{4}+2e_{5}\equiv 
			\begin{tabular}{|l|}
				\hline
				$0$ \\ 
				$1$ \\ 
				$4$ \\ 
				$4$ \\ 
				$1$ \\ 
				$2$%
			\end{tabular}%
		\end{eqnarray*}
		
		\item If $u\equiv 
		\begin{tabular}{|l|}
		\hline
		$+2$ \\ 
		$-1$ \\ 
		$+3$%
		\end{tabular}%
		,$ then:%
		\begin{equation*}
		u\equiv 2e_{0}-e_{1}+3e_{2}\equiv 2\alpha ^{0}-\alpha ^{1}+3\alpha
		^{2}=2-\alpha +3\alpha ^{2}
		\end{equation*}
	\end{itemize}
\end{example}
\begin{definition}
	A polynumber sequence $\left[ u\left( k\right) \right\rangle _{i}^{\ell }$
	is defined by a polynumber $u$, an initial value $0\leq i\leq \ell ,$ and a
	final value $0\leq f\leq \ell $ such that:%
	\begin{equation*}
	\left[ u\left( k\right) \right\rangle _{i}^{f}\equiv
	\sum\limits_{k=i}^{f}u_{k}e_{k}
	\end{equation*}
\end{definition}

\begin{example}
	The polynomial sequences $\left[ u\left( k\right) \right\rangle _{2}^{4}=%
	\left[ 2-k+3k^{2}\right\rangle _{2}^{4}$ is congruent to%
	\begin{equation*}
	\left[ u\left( k\right) \right\rangle _{1}^{4}\equiv \text{ }u\left(
	2\right) e_{2}+u\left( 3\right) e_{3},\text{ }u\left( 4\right) e_{4}\equiv %
	\left[ 12,\text{ }26,\text{ }46\right\rangle 
	\end{equation*}%
	\text{ \ \ }We want to consider the finite sequences that keep going on, in
	a limited way (finite computable way).\newline
\end{example}

\begin{remark}
	From this definition we get the following congruences:
	
	\begin{itemize}
		\item $\left[ \left[ 0\right] \right] \equiv e_{0}=\left[ 
		\begin{array}{cccc}
		1 & 0 & 0 & \cdots%
		\end{array}%
		\right\rangle \equiv \alpha ^{0}=%
		\begin{tabular}{|l|}
		\hline
		$1$ \\ 
		$0$ \\ 
		$\vdots $%
		\end{tabular}%
		$
		
		\item $\left[ \left[ 1\right] \right] \equiv e_{1}=\left[ 
		\begin{array}{cccc}
		0 & 1 & 0 & \cdots%
		\end{array}%
		\right\rangle \equiv \alpha ^{1}=%
		\begin{tabular}{|l|}
		\hline
		$0$ \\ 
		$1$ \\ 
		$\vdots $%
		\end{tabular}%
		\equiv \alpha $
		
		\item $\left[ \left[ k\right] \right] \equiv e_{k}=\left[ 
		\begin{array}{cccccc}
		0 & \cdots & 1 & \cdots &  & 
		\end{array}%
		\right\rangle \equiv \alpha ^{k}=$ $%
		\begin{tabular}{|c|}
		\hline
		$0$ \\ 
		$\vdots $ \\ 
		$1$ \\ 
		$\vdots $%
		\end{tabular}%
		$
		
		\item $u=\left[ 
		\begin{array}{ccccccccc}
		u_{0} & u_{1} & \cdots  & u_{k-1} & u_{k} & u_{k+1} & \cdots  & u_{\ell -1}
		& u_{\ell }%
		\end{array}%
		\right\rangle \equiv 
		\begin{tabular}{|c|}
		\hline
		$u_{0}$ \\ 
		$u_{1}$ \\ 
		$\vdots $ \\ 
		$u_{k-1}$ \\ 
		$u_{k}$ \\ 
		$u_{k+1}$ \\ 
		$\vdots $ \\ 
		$u_{\ell -1}$ \\ 
		$u_{\ell }$%
		\end{tabular}%
		$
	\end{itemize}
\end{remark}

\begin{definition}
	A polyserie $u\left( \alpha \right) $ (or $\alpha $ powet serie) is the sum
	of $\alpha $ powers of \ usual numbers, such that:%
	\begin{equation*}
	u\left( \alpha \right) =\left[ u_{k}\alpha ^{k}\right\rangle _{k\in 
		\mathbb{N}
	}=u_{0}+u_{1}\alpha +u_{2}\alpha ^{2}+\cdots \equiv 
	\begin{tabular}{|c|}
	\hline
	$u_{0}$ \\ 
	$u_{1}$ \\ 
	$u_{2}$ \\ 
	$\vdots $%
	\end{tabular}%
	\end{equation*}%
	endowed with the addition operation:%
	\begin{equation*}
	\text{ \ \ }u+v=\left[ u_{k}+v_{k}\right\rangle _{k\in 
		\mathbb{N}
	}\equiv 
	\begin{tabular}{|c|}
	\hline
	$u_{0}+v_{0}$ \\ 
	$u_{1}+v_{1}$ \\ 
	$u_{2}+v_{2}$ \\ 
	$\vdots $%
	\end{tabular}%
	\end{equation*}%
	and the multiplication operation:%
	\begin{equation*}
	\text{\ \ }u\text{ }.v=\left[ \sum\limits_{i+j=k}u_{i}v_{j}\right\rangle
	_{k\in 
		\mathbb{N}
	}\equiv 
	\begin{tabular}{|c|}
	\hline
	$u_{0}v_{0}$ \\ 
	$u_{1}v_{0}+u_{1}v_{0}$ \\ 
	$u_{2}v_{0}+u_{1}v_{1}+v_{2}u_{0}$ \\ 
	$\vdots $%
	\end{tabular}%
	\end{equation*}
\end{definition}
\begin{example}
	\text{ \ \ }
	
	\begin{itemize}
		\item The ongoing polyserie $\left[ 2-k+3k^{2}\right\rangle _{1}^{+}$ is a
		clip of an ongoing sequence for the polynumber $2-k+3k^{2}$ starting from
		the integer $1$%
		\begin{equation*}
		\left[ 2-k+3k^{2}\right\rangle _{1}^{+}\equiv \left[ 4,\text{ }12,\text{ }26,%
		\text{ }46,\cdots \right\rangle 
		\end{equation*}%
		\textbf{See }\cite{njwildbergerAlgCal}.
		
		\item We have $\ $%
		\begin{equation*}
		\begin{tabular}{|c|}
		\hline
		$0$ \\ 
		$1$ \\ 
		$2$ \\ 
		$3$ \\ 
		$\vdots $%
		\end{tabular}%
		+%
		\begin{tabular}{|c|}
		\hline
		$1$ \\ 
		$2$ \\ 
		$4$ \\ 
		$8$ \\ 
		$\vdots $%
		\end{tabular}%
		\equiv \left( 0+1.\alpha +2\alpha ^{2}+3\alpha ^{3}+\cdots \right) +\left(
		1+2.\alpha +2\alpha ^{2}+8\alpha ^{3}\cdots \right) =\left( 1+3.\alpha
		+6\alpha ^{2}+11\alpha ^{3}+\cdots \right) \equiv 
		\begin{tabular}{|c|}
		\hline
		$1$ \\ 
		$3$ \\ 
		$6$ \\ 
		$11$ \\ 
		$\vdots $%
		\end{tabular}%
		\end{equation*}%
		and 
		\begin{equation*}
		\begin{tabular}{|c|}
		\hline
		$1$ \\ 
		$1$ \\ 
		$1$ \\ 
		$1$ \\ 
		$1$ \\ 
		$\vdots $%
		\end{tabular}%
		\cdot 
		\begin{tabular}{|c|}
		\hline
		$0$ \\ 
		$1$ \\ 
		$3$ \\ 
		$7$ \\ 
		$15$ \\ 
		$\vdots $%
		\end{tabular}%
		\equiv \left( 1+\alpha +\alpha ^{2}+\alpha ^{3}+\alpha ^{4}+\cdots \right)
		\cdot \left( \alpha +3\alpha ^{2}+7\alpha ^{3}+15\alpha ^{4}+\cdots \right)
		=\alpha +4\alpha ^{2}+11\alpha ^{3}+26\alpha ^{4}+\cdots \equiv 
		\begin{tabular}{|c|}
		\hline
		$0$ \\ 
		$1$ \\ 
		$4$ \\ 
		$11$ \\ 
		$26$ \\ 
		$\vdots $%
		\end{tabular}%
		\end{equation*}
	\end{itemize}
\end{example}
\begin{algorithm}[Square root]
	\text{ \ \ }\newline
	We search a solution for $\beta $ of the equation:%
	\begin{equation*}
	\beta ^{2}\equiv 1+\alpha ,
	\end{equation*}%
	such that:%
	\begin{equation*}
	\alpha =\left[ 
	\begin{array}{cccc}
	0 & 1 & 0 & \cdots%
	\end{array}%
	\right\rangle .
	\end{equation*}%
	We have:%
	\begin{equation*}
	\begin{array}{ccl}
	\beta ^{2}\equiv 1+\alpha & \Leftrightarrow & \left[ 
	\begin{array}{ccccc}
	\beta _{0} & \beta _{1} & \beta _{2} & \beta _{3} & \cdots%
	\end{array}%
	\right\rangle \left[ 
	\begin{array}{ccccc}
	\beta _{0} & \beta _{1} & \beta _{2} & \beta _{3} & \cdots%
	\end{array}%
	\right\rangle \equiv 1+\alpha \\ 
	& \Leftrightarrow & \left[ 
	\begin{array}{cccc}
	\beta _{0}^{2} & 2\beta _{0}\beta _{1} & \beta _{1}^{2}+2\beta _{0}\beta _{2}
	& \cdots%
	\end{array}%
	\right\rangle \equiv 1+\alpha .%
	\end{array}%
	\end{equation*}%
	We deduce that: 
	\begin{equation*}
	\left\{ 
	\begin{array}{l}
	\beta _{0}^{2}\equiv 1 \\ 
	2\beta _{0}\beta _{1}\equiv 1 \\ 
	\beta _{1}^{2}+2\beta _{0}\beta _{2}\equiv 0 \\ 
	\beta _{0}\beta _{3}+\beta _{1}\beta _{2}+\beta _{2}\beta _{1}+\beta
	_{3}\beta _{0}\equiv 0%
	\end{array}%
	\right.
	\end{equation*}%
	Then the Solution of $\beta $ is given by its coordinates:%
	\begin{eqnarray*}
		&&\left[ \beta _{0}\equiv -1,\text{ }\beta _{1}\equiv -\frac{1}{2},\text{ }%
		\beta _{2}\equiv \displaystyle\frac{1}{8},\text{ }\beta _{3}\equiv -\frac{1}{%
			16}\right] , \\
		&&\text{or} \\
		&&\left[ \beta _{0}\equiv 1,\text{ }\beta _{1}\equiv \displaystyle\frac{1}{2}%
		,\text{ }\beta _{2}\equiv -\displaystyle\frac{1}{8},\text{ }\beta _{3}\equiv %
		\displaystyle\frac{1}{16}\right]
	\end{eqnarray*}
\end{algorithm}
\begin{definition}[Informally]
	\text{ \ \ }\newline
	A sequence is a set of related events, movements, or items that follow each
	other in a particular order.
\end{definition}

\begin{example}
	\text{ \ \ }
	
	\begin{itemize}
		\item Integer sequences: 
		\begin{equation*}
		S=\left( a_{i}\right) _{i\in 
			\mathbb{N}
		}=\left( a_{0},\text{ }a_{1},\text{ }a_{2},\ldots \right)
		\end{equation*}
		
		\item Forward difference operator of a sequence: 
		\begin{equation*}
		\Delta \left( S\right) =\left( a_{i+1}-a_{i}\right) _{i\in 
			\mathbb{N}
		}=\left( a_{1}-a_{0},\text{ }a_{2}-a_{1},\ldots \right)
		\end{equation*}
		
		\item Summation operator of a sequence: 
		\begin{equation*}
		\textstyle\sum \left( S\right) =\left( \textstyle\sum_{k<i}a\right) _{i\in 
			\mathbb{N}
		}=\left( 0,\text{ }a_{0},\text{ }a_{0}+a_{1},\text{ }a_{0}+a_{1}+a_{2},%
		\ldots \right)
		\end{equation*}
	\end{itemize}
	
	Thomas Harriot (1560-1621) was the first to introduce the difference tables: 
	\begin{table}[H]
		\caption{Harriot's polyseries difference tables.}\centering \vspace{0.5cm}
		\begin{tabular}{|c|c|c|c|c|c|}
			\hline
			$S=\langle n^{4}\rangle _{n\in 
				\mathbb{N}
				_{\geq 0}}$ & $\Delta (S)$ & $\Delta ^{2}(S)$ & $\Delta ^{3}(S)$ & $\Delta
			^{4}(S)$ & $\Delta ^{5}(S)$ \\ \hline
			$0$ & $1$ & $14$ & $36$ & $24$ & $0$ \\ \hline
			$1$ & $15$ & $50$ & $60$ & $24$ & $0$ \\ \hline
			$16$ & $65$ & $110$ & $84$ & $24$ & $\vdots $ \\ \hline
			$81$ & $175$ & $194$ & $108$ & $\vdots $ & $\vdots $ \\ \hline
			$256$ & $369$ & $302$ & $\vdots $ & $\vdots $ & $\vdots $ \\ \hline
			$625$ & $671$ & $\vdots $ & $\vdots $ & $\vdots $ & $\vdots $ \\ \hline
			$1296$ & $\vdots $ & $\vdots $ & $\vdots $ & $\vdots $ & $\vdots $ \\ \hline
			$\vdots $ & $\vdots $ & $\vdots $ & $\vdots $ & $\vdots $ & $\vdots $ \\ 
			\hline
		\end{tabular}%
	\end{table}
\end{example}
\normalsize
\subsection{Harriot's triangular summation tables}

\begin{itemize}
	\item The columns of \ the ongoing Harriot's polyseries summation table are: 
	\begin{table}[H]
		\caption{Harriot's Repeated polyseries summations of $S=\langle 1\rangle $.}%
		\vspace{0.5cm}\centering 
		\begin{tabular}{|c|c|c|c|c|c|c|}
			\hline
			$S=\langle 1\rangle _{n\in 
				\mathbb{N}
				_{\geq 0}}$ & $\sum (S)$ & $\sum^{2}(S)$ & $\sum^{3}(S)$ & $\sum^{4}(S)$ & $%
			\sum^{5}(S)$ &  \\ \hline
			$1$ & $0$ & $0$ & $0$ & $0$ & $0$ & $\cdots $ \\ \hline
			$1$ & $1$ & $0$ & $0$ & $0$ & $0$ & $\cdots $ \\ \hline
			$1$ & $2$ & $1$ & $0$ & $0$ & $0$ & $\cdots $ \\ \hline
			$1$ & $3$ & $3$ & $1$ & $0$ & $0$ & $\cdots $ \\ \hline
			$1$ & $4$ & $6$ & $4$ & $1$ & $0$ & $\cdots $ \\ \hline
			$1$ & $5$ & $10$ & $10$ & $5$ & $1$ & $\cdots $ \\ \hline
			$1$ & $6$ & $15$ & $20$ & $15$ & $6$ & $\cdots $ \\ \hline
			$\vdots $ & $7$ & $21$ & $35$ & $35$ & $21$ & $\cdots $ \\ \hline
			$\vdots $ & $8$ & $28$ & $56$ & $70$ & $56$ & $\cdots $ \\ \hline
			$\vdots $ & $\vdots $ & $\vdots $ & $\vdots $ & $\vdots $ & $\vdots $ & $%
			\ddots $ \\ \hline
		\end{tabular}%
	\end{table}
	\normalsize
	\item The initial ongoing Harriot's polyserie $h_{0}\equiv \left[
	1\right\rangle =1,\text{ }1,\text{ }1,\cdots $
	
	\item The first summation of it is $h_{1}\equiv \left[ n\right\rangle _{n\in 
		\mathbb{N}
	}=0,\text{ }1,\text{ }2,\cdots $
	
	\item Its second one is%
	\begin{equation*}
	h_{2}\equiv \left[ \binom{n}{2}\right\rangle _{n\in 
		\mathbb{N}
	}=\left[ \frac{n\left( n-1\right) }{2!}\right\rangle _{n\in 
		\mathbb{N}
	}=0,\text{ }0,\text{ }1,\text{ }3,\text{ }6,\text{ }10,\cdots \ \ \ \ \ 
	\text{(Triangular Numbers)}.
	\end{equation*}
	
	\item The third is:%
	\begin{equation*}
	h_{3}\equiv \left[ \binom{n}{3}\right\rangle _{n\in 
		\mathbb{N}
	}=\left[ \frac{n\left( n-1\right) \left( n-2\right) }{3!}\right\rangle
	_{n\in 
		\mathbb{N}
	}=0,\text{ }0,\text{ }0,\text{ }1,\text{ }4,\text{ }10,\cdots \text{%
		(Pyramidal Numbers)}.
	\end{equation*}
	
	\item The general term of the ongoing Harriot's polyserie is: 
	\begin{equation*}
	h_{n,k}=\left[ \binom{n}{k}\right\rangle _{n\in 
		\mathbb{N}
	}=\left[ \frac{n\left( n-1\right) \cdots \left( n-k+1\right) }{k!}%
	\right\rangle _{n\in 
		\mathbb{N}
	}=\left[ \frac{n^{\underline{k}}}{k!}\right\rangle _{n\in 
		\mathbb{N}
	}
	\end{equation*}%
	{\ $n^{\underline{k}}=n\left( n-1\right) \cdots \left( n-k+1\right) $ is the 
		$n$ to the $k$ falling (Knuth. D. notation) }
	
	\item Harriot's Triangular (Binomial) sequences $h_{0},h_{1},h_{2}\cdots $,
	by using the binomial formula: 
	\begin{equation*}
	h_{n,k}=h_{k}\left( n\right) =\binom{n}{k},\text{ \ \ }k,\text{ }n=0,\text{ }%
	1,\text{ }2,\cdots
	\end{equation*}%
	We get the following table:
	\begin{table}[H]
		\caption{Array of $h_{n,k}$ (rows indexed by $n$, columns by $k$).}\vspace{0.5cm}\centering%
		\begin{tabular}{c|p{1.2cm}p{1.2cm}p{1.2cm}p{1.2cm}p{1.2cm}p{1.2cm}}
			& \multicolumn{6}{c}{$k$} \\ \cline{2-7}
			$n$ & $1$ & $2$ & $3$ & $4$ & $5$ & $6$ \\ \hline
			& $1$ & 0 & 0 & 0 & $\cdots $ & $\cdots $ \\ \hline
			& $1$ & 1 & 0 & 0 & $\cdots $ & $\cdots $ \\ \hline
			& $1$ & 2 & 1 & 0 & $\cdots $ & $\cdots $ \\ \hline
			& $1$ & 3 & 3 & 1 & $\cdots $ & $\cdots $ \\ \hline
			& $1$ & 4 & 6 & 4 & 1 & $\cdots $ \\ \hline
			$\vdots $ & $\vdots $ & $\vdots $ & $\vdots $ & $\vdots $ & $\vdots $ & $%
			\ddots $ \\ \hline
		\end{tabular}%
	\end{table}
\end{itemize}
\normalsize

\begin{theorem}[Harriot's Difference Theorem]
	\begin{equation*}
	\Delta \left( h_{k}\right) =h_{k-1}\text{ \ \ \ \ \ }\left( \text{For \ }%
	k\in 
	\mathbb{N}
	_{\geq 1}\right) \text{\ }
	\end{equation*}
\end{theorem}

\begin{proof}
	\text{ \ \ }\newline
	By definition:%
	\begin{equation*}
	h_{k}\left( n\right) =\binom{n}{k}.
	\end{equation*}
	So we get: 
	\begin{equation*}
	\Delta \left( h_{k}\right) \left( n\right) =h_{k}\left( n+1\right)
	-h_{k}\left( n\right) =\binom{n+1}{k}-\binom{n}{k}=\binom{n}{k-1}%
	=h_{k-1}\left( n\right)
	\end{equation*}
\end{proof}

\begin{corollary}
	\text{ \ \ }\newline
	For any sequence $S=a_{0},a_{1},a_{2}\cdots $ we can generate difference or
	summation table.
\end{corollary}

\begin{example}
	\text{ \ \ }\newline
	From the sequence:%
	\begin{equation*}
	S=1,\text{ }1,\text{ }3,\text{ }13,\text{ }37,\text{ }81,\text{ }212,
	\end{equation*}%
	we generate the following summation table: 
	\begin{table}[H]
		\caption{Upper triangular integer matrix }\vspace{0.5cm}\centering
		\begin{tabular}{ccccccccccc}
			\cline{1-9}
			\multicolumn{1}{|c}{0} & \multicolumn{1}{|c}{0} & \multicolumn{1}{|c}{6} & 
			\multicolumn{1}{|c}{2} & \multicolumn{1}{|c}{0} & \multicolumn{1}{|c}{1} & 
			\multicolumn{1}{|c}{0} & \multicolumn{1}{|c}{0} & \multicolumn{1}{|c}{0} & 
			\multicolumn{1}{|c}{} &  \\ \cline{1-9}
			\multicolumn{1}{|c}{0} & \multicolumn{1}{|c}{0} & \multicolumn{1}{|c}{6} & 
			\multicolumn{1}{|c}{8} & \multicolumn{1}{|c}{2} & \multicolumn{1}{|c}{1} & 
			\multicolumn{1}{|c}{1} & \multicolumn{1}{|c}{0} & \multicolumn{1}{|c}{0} & 
			\multicolumn{1}{|c}{} &  \\ \cline{1-9}
			\multicolumn{1}{|c}{0} & \multicolumn{1}{|c}{0} & \multicolumn{1}{|c}{6} & 
			\multicolumn{1}{|c}{14} & \multicolumn{1}{|c}{10} & \multicolumn{1}{|c}{3} & 
			\multicolumn{1}{|c}{2} & \multicolumn{1}{|c}{1} & \multicolumn{1}{|c}{0} & 
			\multicolumn{1}{|c}{} &  \\ \cline{1-9}
			\multicolumn{1}{|c}{0} & \multicolumn{1}{|c}{0} & \multicolumn{1}{|c}{0} & 
			\multicolumn{1}{|c}{20} & \multicolumn{1}{|c}{24} & \multicolumn{1}{|c}{13}
			& \multicolumn{1}{|c}{5} & \multicolumn{1}{|c}{3} & \multicolumn{1}{|c}{1} & 
			\multicolumn{1}{|c}{} &  \\ \cline{1-9}
			\multicolumn{1}{|c}{0} & \multicolumn{1}{|c}{0} & \multicolumn{1}{|c}{0} & 
			\multicolumn{1}{|c}{0} & \multicolumn{1}{|c}{44} & \multicolumn{1}{|c}{37} & 
			\multicolumn{1}{|c}{18} & \multicolumn{1}{|c}{8} & \multicolumn{1}{|c}{4} & 
			\multicolumn{1}{|c}{} &  \\ \cline{1-9}
			\multicolumn{1}{|c}{0} & \multicolumn{1}{|c}{0} & \multicolumn{1}{|c}{0} & 
			\multicolumn{1}{|c}{0} & \multicolumn{1}{|c}{0} & \multicolumn{1}{|c}{81} & 
			\multicolumn{1}{|c}{55} & \multicolumn{1}{|c}{26} & \multicolumn{1}{|c}{12}
			& \multicolumn{1}{|c}{} &  \\ \cline{1-9}
			\multicolumn{1}{|c}{0} & \multicolumn{1}{|c}{0} & \multicolumn{1}{|c}{0} & 
			\multicolumn{1}{|c}{0} & \multicolumn{1}{|c}{0} & \multicolumn{1}{|c}{0} & 
			\multicolumn{1}{|c}{131} & \multicolumn{1}{|c}{81} & \multicolumn{1}{|c}{38}
			& \multicolumn{1}{|c}{} &  \\ \cline{1-9}\cline{1-9}
			\multicolumn{1}{|c}{0} & \multicolumn{1}{|c}{0} & \multicolumn{1}{|c}{$0$} & 
			\multicolumn{1}{|c}{0} & \multicolumn{1}{|c}{0} & \multicolumn{1}{|c}{0} & 
			\multicolumn{1}{|c}{0} & \multicolumn{1}{|c}{212} & \multicolumn{1}{|c}{219$%
				\;\;$} & \multicolumn{1}{|c}{} & $\underrightarrow{\sum }$ \\ \cline{1-9}
			&  & $\underleftarrow{\Delta }$ &  &  &  &  &  &  &  & 
		\end{tabular}%
	\end{table}
\end{example}

\normalsize
\subsection{Chains}

\begin{definition}
	\text{ \ \ }\newline
	For a fixed counting number $k$ and an integer $n$, the $k$-chain on \ $n$
	is 
	\begin{equation*}
	C\left( k,n\right) \equiv n,\text{ }n+1,\text{ }\cdots ,\text{ }n+k-1
	\end{equation*}
\end{definition}

\begin{example}
	{\ \text{} }
	
	\begin{itemize}
		\item The $3$-chain on $17$ is $17,18,19$
		
		\item The $1$-chain on $-51$ is $-51$
		
		\item The $8$-chain on $-2$ is $-2,$ $-1,$ $0,$ $1,$ $2,$ $3,$ $4,$ $5$
	\end{itemize}
\end{example}

\subsection{Integral sequences}

\begin{definition}
	\text{ \ \ }\newline
	An integral $k-$sequence is an explicit assignment of integers to the
	elements of a $k$-chain
\end{definition}

\begin{example}
	\text{ \ \ }\newline
	From this table: 
	\begin{table}[H]
		\caption{Values of $S(n)$ for integer shifts}\vspace{0.5cm}\centering
		\begin{tabular}{|c|c|c|c|c|c|c|}
			\hline
			$n$ & $-1$ & $0$ & $+1$ & $+2$ & $+3$ & $+4$ \\ \hline
			$S(n)$ & $+1$ & $+2$ & $-1$ & $+3$ & $0$ & $+2$ \\ \hline
		\end{tabular}%
	\end{table}
	\normalsize
	We deduce the sequence:%
	\begin{eqnarray*}
		S &=&1_{-1},\text{\ \ }2_{0},\text{ }-1_{1},\text{ }3_{2},\text{ }0_{3},%
		\text{ }2_{4} \\
		S_{-1}^{4} &=&1,\text{ }2,\text{ }-1,\text{ }3,\text{ }0,\text{ }2\text{ \ \
			the limits are: }-1\text{ and }4\text{, \ \ \ the size is }6
	\end{eqnarray*}
\end{example}

\begin{remark}
	\text{ \ \ }
	
	\begin{itemize}
		\item {\ The default beginning limit is $0.$ }
		
		\item {\ The term sequence means $k$-sequence for some counting number $k.$ }
		
		\item {\ We only consider finite sequences. }
	\end{itemize}
\end{remark}

\subsection{Clips}

\begin{definition}
	\text{ \ \ }\newline
	A clip is a representation of a (consecutive) portion of a sequence
\end{definition}

\begin{example}
	The sequence $S=1_{-1},$ $2_{0},$ $-1_{1},$ $3_{2},$ $0_{3},$ $2_{4}$ has
	clips:
	
	\begin{itemize}
		\item {\ $\ell =1_{-1},$ $2_{0},$ $-1_{1},$ $3_{2},\cdots $(left clip) }
		
		\item {\ $j=\cdots ,3_{2},$ $0_{3},$ $2_{4}$ \ \ (right clip) }
		
		\item {\ $k=\cdots ,-1,$ $3_{2},\cdots $ (double clip) }
		
		\item {\ $s=1_{-1},$ $2_{0},$ $-1_{1},$ $3_{2},$ $0_{3},$ $2_{4}$ (total
			clip) }
	\end{itemize}
\end{example}
\begin{remark}{\ \text{} }
	\begin{itemize}
		\item {\ A clip must have at least one element. }
		\item {\ The default clip kind is the left one starting at 0. }
		\item {\ A given sequence will have many clips: partial representation of
			the sequence. }
	\end{itemize}
\end{remark}
\begin{example}
	From the sequence: $S=10,\text{ }9,\text{ }8,\text{ }7,\text{ }6,\text{ }5,%
	\text{ }4, $ we get: 
	\begin{table}[H]
		\caption{ Clip and assignment data classification}\vspace{0.5cm}\centering
		\begin{tabular}{|@{}c|cc|}
			\hline
			\textbf{clip} &  & \textbf{assignment data} \\ \hline
			\multicolumn{1}{|c|}{$10,\,9,\,8,\ldots $} & \multicolumn{1}{|c}{} & 
			\multicolumn{1}{c|}{right data (accepted)} \\ \hline
			\multicolumn{1}{|c|}{$\ldots ,7,\,6,\ldots $} & \multicolumn{1}{|c}{} & 
			\multicolumn{1}{c|}{wrong data (non accepted)} \\ \hline
			\multicolumn{1}{|c|}{$\ldots ,\,7_{3},\,6_{4},\ldots $} & 
			\multicolumn{1}{|c}{} & \multicolumn{1}{c|}{right data} \\ \hline
			\multicolumn{1}{|c|}{$\ldots ,\,4_{6}$} & \multicolumn{1}{|c}{} & 
			\multicolumn{1}{c|}{right} \\ \hline
			\multicolumn{1}{|c|}{$\ldots $} & \multicolumn{1}{|c}{} & 
			\multicolumn{1}{c|}{wrong} \\ \hline
		\end{tabular}%
	\end{table}
\end{example}
\normalsize

\subsection{Binomial Chu-Vandermonde Identity}

\begin{definition}
We define the powers of a polynumber $u$ for a positive natural number $n\geq 1$ by: 
\begin{eqnarray*}
	u^{0:\text{ }n} &\coloneqq&1, \\
	u^{n\text{ }:\text{ }0} &\coloneqq&u^{n}, \\
	u^{n\text{ }:\text{ }1} &\coloneqq&u^{\overline{n}}=\prod_{k=0}^{n-1}\left(
	u+k\right) =u\left( u+1\right) \left( u+2\right) \cdots \left( u+\left(
	n-1\right) \right),  \\
	u^{n\text{ }:\text{ }-1} &\coloneqq&u^{\underline{n}}=\prod_{k=0}^{n-1}\left(
	u-k\right) =u\left( u-1\right) \left( u-2\right) \cdots \left( u-\left(
	n-1\right) \right), 
\end{eqnarray*}

\bigskip and we define the Ladder powers of a polynumber $u$  for a polynumber $t$ by:%
\begin{eqnarray*}
	u^{0\text{ }:\text{ }t} &\coloneqq&1,\text{ } \\
	u^{1\text{ }:\text{ }t} &\coloneqq&u,\text{ } \\
	u^{2\text{ }:\text{ }t} &\coloneqq&u\left( u+t\right) , \\
	&&\vdots  \\
	u^{n\text{ }:\text{ }t} &\coloneqq&\prod_{k=0}^{n-1}\left( u+kt\right) =u\left(
	u+t\right) \left( u+2t\right) \cdots \left( u+\left( n-1\right) t\right). 
\end{eqnarray*}
\end{definition}

\begin{theorem}[Binomial Chu-Vandermonde General identity]
	\text{ \ \ }\newline 
	For any natural number $n$ and polynumbers (eventually rationals as special particular case) $u,v,$ and $t$%
	\begin{equation}
	\left( u+v\right) ^{n\text{ }:\text{ }t}=\sum\limits_{k=0}^{n}\binom{n}{k}%
	u^{\left( n-k\right) \text{ }:\text{ }t}\, . \,v^{k\text{ }:\text{ }t}\label{Chu-Vandermonde}
	\end{equation}
\end{theorem}

	\begin{proof}
	
		For any natural number \(n\) and polynumbers \(u,v\) and step \(t\), we prove the identity by induction on \(n\).

	 \setlength{\LTleft}{-0.5in} 
	  \setlength{\LTright}{0pt} 
	\begin{prooftable}
		1 & Base case \(n=0\): \((u+v)^{0:t}=1\) and the sum is \(\binom{0}{0}u^{0:t}v^{0:t}=1\). & Definition of \(x^{0:t}\) and binomial coefficient \(\binom{0}{0}=1\). \\
		2 & Induction hypothesis: assume the identity holds for some \(n\ge0\). & Standard induction setup. \\
		3 & Consider \( (u+v)^{(n+1):t} = (u+v)^{n:t}\cdot\bigl(u+v + n t\bigr). \) & By definition \(x^{(n+1):t}=x^{n:t}\,(x+n t)\) with \(x=u+v\). \\
		4 & Substitute the induction hypothesis for \((u+v)^{n:t}\): \\
		& \(\displaystyle (u+v)^{(n+1):t}
		=\Biggl(\sum_{k=0}^{n}\binom{n}{k}u^{(n-k):t}v^{k:t}\Biggr)\cdot\bigl(u+v + n t\bigr).\) & Use step 2. \\
		5 & Observe the linear decomposition \(\;u+v+n t=(u+(n-k)t)+(v+k t)\;\) for each \(k\). & Simple algebra: add and subtract \(k t\). \\
		6 & Multiply termwise and split the sum into two sums:

$(u+v)^{(n+1):t}=\sum_{k=0}^{n}\binom{n}{k}u^{(n-k):t}\bigl(u+(n-k)t\bigr)\,v^{k:t}+\sum_{k=0}^{n}\binom{n}{k}u^{(n-k):t}\,v^{k:t}\bigl(v+k t\bigr).$
		
		& Distribute \((u+v+n t)\) using step 5. \\
		7 & Use the defining recurrence \(x^{m+1:t}=x^{m:t}\,(x+(m-1)t)\) (equivalently \(u^{(n-k+1):t}=u^{(n-k):t}\,(u+(n-k)t)\)) to rewrite the first sum:

		$
		\sum_{k=0}^{n}\binom{n}{k}u^{(n-k):t}\bigl(u+(n-k)t\bigr)\,v^{k:t}
		=\sum_{k=0}^{n}\binom{n}{k}u^{(n-k+1):t}\,v^{k:t}.
	$
		
		& Apply the step‑extension formula for \(u^{\cdot :t}\). \\
		8 & Similarly rewrite the second sum using \(v^{k+1:t}=v^{k:t}\,(v+k t)\):

		$
		\sum_{k=0}^{n}\binom{n}{k}u^{(n-k):t}\,v^{k:t}\bigl(v+k t\bigr)
		=\sum_{k=0}^{n}\binom{n}{k}u^{(n-k):t}\,v^{k+1:t}.
		$
		
		& Apply the step‑extension formula for \(v^{\cdot :t}\). \\
		9 & Reindex the second sum by \(j=k+1\). Then

		$
		\sum_{k=0}^{n}\binom{n}{k}u^{(n-k):t}\,v^{k+1:t}
		=\sum_{j=1}^{n+1}\binom{n}{j-1}u^{(n-(j-1)):t}\,v^{j:t}.
		$
		
		& Change of index \(j=k+1\). \\
		10 & Combine the two sums (from steps 7 and 9):

		$
		(u+v)^{(n+1):t}
		=\sum_{k=0}^{n}\binom{n}{k}u^{(n-k+1):t}\,v^{k:t}
		+\sum_{k=1}^{n+1}\binom{n}{k-1}u^{(n-k+1):t}\,v^{k:t}.
		$
		
		& Align indices so both sums are over the same power pattern. \\
		11 & Merge the two sums termwise using the binomial recurrence \(\binom{n}{k}+\binom{n}{k-1}=\binom{n+1}{k}\):
	$(u+v)^{(n+1):t}=\sum_{k=0}^{n+1}\binom{n+1}{k}u^{(n+1-k):t}\,v^{k:t}.
		$
		
		& Boundary terms: for \(k=0\) the second sum contributes nothing; for \(k=n+1\) the first sum contributes nothing; interior terms combine by the binomial identity. \\
		12 & This completes the induction, so the identity holds for all \(n\in\mathbb{N}\). & Induction principle. \\
	\end{prooftable}
	\end{proof}
\normalsize
\begin{example}[Newton polyseries]
	\text{ \ \ }\newline
As a particular case from \eqref{Chu-Vandermonde} Chu-Vandermonde general identity, for any two rationals $r,$ $s,$ and a polynumber $u$ with $t=-1,$ we get:%
\begin{equation}
\left( \sum\limits_{k=0}^{+}\frac{r^{k\text{ }:\text{ }-1}}{k!}u^{k}\right)
\,.\,\left( \sum\limits_{\ell =0}^{+}\frac{s^{\ell \text{ }:\text{ }-1}}{\ell
	!}u^{\ell }\right) =\left( \sum\limits_{n =0}^{+}\frac{\left( r+s\right)
	^{n\text{ }:\text{ }-1}}{\left( r+s\right) !}u^{n}\right),\label{powersum}
\end{equation}%
Such that:
	\begin{equation*}
	\left( r+s\right) ^{n\text{ }:\text{ }-1}=\sum\limits_{k=0}^{n}\binom{n}{k}r^{\left( n-k\right) \text{ }:\text{ }%
		-1}\,.\,s^{k\text{ }:\text{ }-1}
	\end{equation*}%
	We define:
	\begin{equation*}
	\left( 1+u\right\rangle ^{r}\coloneqq \sum\limits_{\ n=0}^{+}\frac{r^{n\text{ }:\text{
			}-1}}{n!}u^{n}\quad \text{and}\quad \left( 1+u\right\rangle
	^{s}\coloneqq \sum\limits_{\ n=0}^{+}\frac{s^{n\text{ }:\text{ }-1}}{n!}u^{n}.
	\end{equation*}%
	
Using \eqref{powersum}, we deduce that: 

\[
\left( 1+u\right\rangle ^{r}\,.\,\left( 1+u\right\rangle ^{s}=\left(
\sum\limits_{\ n=0}^{+}\frac{r^{n\text{ }:\text{ }-1}}{n!}u^{n}\right)
\,.\,\left( \sum\limits_{\ n=0}^{+}\frac{s^{n\text{ }:\text{ }-1}}{n!}%
u^{n}\right) =\sum\limits_{\ n=0}^{+}\left( \sum\limits_{\ k=0}^{n}\frac{r^{k%
		\text{ }:\text{ }-1}}{k!}\frac{s^{n-k\text{ }:\text{ }-1}}{\left( n-k\right)
	!}\right) u^{n}
\]
\[
=\sum\limits_{\ n=0}^{+}\frac{1}{n!}\left( \sum\limits_{\ k=0}^{n}n!\frac{%
	r^{k\text{ }:\text{ }-1}}{k!}\frac{s^{n-k\text{ }:\text{ }-1}}{\left(
	n-k\right) !}\right) u^{n}=\sum\limits_{\ n=0}^{+}\frac{1}{n!}\left(
\sum\limits_{\ k=0}^{n}\binom{n}{k}r^{k\text{ }:\text{ }-1}.s^{n-k\text{ }:%
	\text{ }-1}\right) u^{n}
\]

\[
=\sum\limits_{\ n=0}^{+}\frac{1}{n!}\left( \left( r+s\right) ^{n:-1}\right)
u^{n}=\left( 1+u\right\rangle ^{r+s}.
\]

\normalsize	
\text{\ }\\

	We get the Newton's Polyseries Identity:
	\begin{equation}
	\forall u\in \mathbf{Polynumbers,}\text{ }r,s\in \mathbf{Rationals,}\text{ \ 
	}\left( 1+u\right\rangle ^{r}\,.\,\left( 1+u\right\rangle ^{s}=\left(
	1+u\right\rangle ^{r+s}  \label{Newton-Identity}
	\end{equation}
\end{example}
\begin{example}[Exponential Euler polyserie]
	\text{ \ \ }\newline
	If $t=0$ we get for any two polynumbers $u,\,v$ the Binomial identity:%
	\begin{equation*}
	\left( u+v\right) ^{n\text{ }:\text{ }0}=\sum\limits_{k=0}^{n}\binom{n}{k}%
	u^{\left( n-k\right) \text{ }:\text{ }0}\,.\,v^{k\text{ }:\text{ }0}.
	\end{equation*}%
For any two rationals $r,$ $s,$ and a polynumber $u,$ we introduce the
exponential polyseries: 
\[
\left( \exp \right\rangle _{u}^{r}=\sum\limits_{n=0}^{+}\frac{r^{n}}{n!}%
u^{n}\quad \text{and}\quad \left( \exp \right\rangle
_{u}^{s}=\sum\limits_{n=0}^{+}\frac{s^{n}}{n!}u^{n}.
\]
Such that: 
\begin{eqnarray*}
	\left( \exp \right\rangle _{u}^{r}\,.\,\left( \exp \right\rangle _{u}^{s}
	&=&\left( \sum\limits_{n=0}^{+}\frac{r^{n\text{ }:\text{ }0}}{n!}%
	u^{n}\right) \left( \sum\limits_{n=0}^{+}\frac{s^{n\text{ }:\text{ }0}}{n!}%
	u^{n}\right)  \\
	&=&\left( \sum\limits_{n=0}^{+}\left( \sum\limits_{\ k=0}^{n}\frac{r^{k:0}}{%
		k!}\frac{s^{n-k\text{ }:\text{ }0}}{\left( n-k\right) !}\right) u^{n}\right) 
	\\
	&=&\left( \sum\limits_{n=0}^{+}\frac{1}{n!}\left( \sum\limits_{\ k=0}^{n}n!%
	\frac{r^{k\text{ }:\text{ }0}}{k!}\frac{s^{n-k\text{ }:\text{ }0}}{\left(
		n-k\right) !}\right) u^{n}\right)  \\
	&=&\left( \sum\limits_{n=0}^{+}\frac{1}{n!}\left( r+s\right)
	^{n:0}u^{n}\right) =\left( \exp \right\rangle _{u}^{r+s}
\end{eqnarray*}

We get the Euler's exponential polyseries identity:%
\begin{equation}
	\forall u\in \mathbf{Polynumbers,}\text{ }r,s\in \mathbf{Rationals,}\text{ \ 
}\left( \exp \right\rangle _{u}^{r}\,.\,\left( \exp \right\rangle
_{u}^{s}=\left( \exp \right\rangle _{u}^{r+s}  \label{Euler-Identity}
\end{equation}
\end{example}
\begin{example}[Newton reciprocal polyserie]

\text{\ }
	We define 
\begin{equation*}
\left( 1-u\right\rangle ^{-r}\coloneqq \sum\limits_{\ n=0}^{+}\frac{r^{n\text{ }:%
		\text{ }1}}{n!}u^{n}
\end{equation*}%
and 
\begin{equation*}
\left( 1-u\right\rangle ^{-s} \coloneqq \sum\limits_{\ n=0}^{+}\frac{s^{n\text{ }:%
		\text{ }1}}{n!}u^{n}.
\end{equation*}%
Then we have:
	\[
\bigskip \left( 1-u\right\rangle ^{-r}\,.\,\left( 1-u\right\rangle
^{-s}=\left( \sum\limits_{\ n=0}^{+}\frac{r^{n\text{ }:\text{ }1}}{n!}%
u^{n}\right) \,.\,\left( \sum\limits_{\ n=0}^{+}\frac{s^{n\text{ }:\text{ }1}%
}{n!}u^{n}\right). 
\]
Since:
	\begin{eqnarray*}
		\left( -s\right) ^{^{n\text{ }:\text{ }1}} &=&\left( -s\right) \left(
		-s+1\right) \left( -s+2\right) \cdots \left( -s+n-1\right) \\
		\ \ &=&\left( -1\right) ^{n}s\left( s-1\right) \left( s-2\right) \cdots
		\left( s-n+1\right) \\
		&=&\left( -1\right) ^{n}s^{n\text{ }:\text{ }-1},
	\end{eqnarray*}%

we get:\scriptsize
	\[ \left( 1-u\right\rangle ^{-r}\,.\,\left( 1-u\right\rangle
	^{-s}=
	\left( \sum\limits_{\ n=0}^{+}\frac{r^{n\text{ }:\text{ -}1}}{n!}\left(
	-1\right) ^{n}u^{n}\right) \,.\,\left( \sum\limits_{\ n=0}^{+}\frac{s^{n%
			\text{ }:\text{ }-1}}{n!}\left( -1\right) ^{n}u^{n}\right) =\sum\limits_{\
		n=0}^{+}\left( \sum\limits_{\ k=0}^{n}\left( -1\right) ^{k}\frac{r^{k\text{ }%
			:\text{ }-1}}{k!}\left( -1\right) ^{n-k}\frac{s^{n-k\text{ }:\text{ }-1}}{%
		\left( n-k\right) !}\right) u^{n}
	\]%
	\[\qquad\qquad\qquad\qquad\quad\quad
	=\sum\limits_{\ n=0}^{+}\left( -1\right) ^{n}\frac{1}{n!}\left( \sum\limits_{\ k=0}^{n}%
	\binom{n}{k}r^{k\text{ }:\text{ }-1}.s^{n-k\text{ }:\text{ }-1}\right)
	u^{n} =\sum\limits_{\ n=0}^{+}\frac{1}{n!}\left( \left( r+s\right)
	^{n:-1}\right) \left( -1\right) ^{n}u^{n}=\sum\limits_{\ n=0}^{+}\frac{1}{n!}%
	\left(\left(- \left( r+s\right)\right)^{n:1}\right) u^{n}.
	\]
\normalsize	
\text{\ }\\
	
	Finally we get:
	 
	\begin{equation}
	\forall u\in \mathbf{Polynumbers,}\text{ }r,s\in \mathbf{Rationals,}\text{ \ 
}	\left( 1-u\right\rangle ^{-r}\,.\,\left( 1-u\right\rangle ^{-s}=\left(
	1-u\right\rangle ^{-\left( r+s\right) }
	\end{equation}
	
\end{example}

\subsection{Arithmetic algorithms and polyseries}

\begin{definition}[Finite Algebra]
	\text{ \ \ }\newline
	A finite algebra (with identity) is a finite vector space $\ \mathbb{A}$
	over the finite modular field $\ \mathbb{F}_{p}$ with a multiplication $a\,
	. \, b$ satisfying for all $a,b,c\in \mathbb{A}$ and $\lambda \in \mathbb{F}%
	_{p}$ the following properties:
	
	\begin{itemize}
		\item \textbf{Associativity} 
		\begin{equation*}
		(a\, . \, b)\, . \, c=a\, . \, (b\, . \, c)
		\end{equation*}
		
		\item \textbf{Distributivity} 
		\begin{equation*}
		\begin{tabular}{l}
		$a\,.(b+c)\;=\;a\,.\,b+a\,.\,c,$ \\ 
		$(a+b)\,.\,c\;=\;a\,.\,c+b\,.\,c,$ \\ 
		$(\lambda \,\ast \,a)\,.\,b=a\,.\,(\lambda \,\ast \,b)\;=\;\lambda \,\ast
		(a\,.\,b)$%
		\end{tabular}%
		\end{equation*}
		
		\item \textbf{Identity} 
		\begin{equation*}
		1\,\,.\, \,a=a\,\,.\,\,1=a
		\end{equation*}
	\end{itemize}
\end{definition}

\begin{definition}
	\text{ \ \ }\newline
	A truncated polyserie is defined finitely as a data-structure (ongoing up to
	a certain finite order $k$)%
	\begin{equation*}
	\text{\ }\alpha _{k}^{0}\equiv \left[ 
	\begin{array}{c}
	1 \\ 
	0 \\ 
	0 \\ 
	\vdots \\ 
	0%
	\end{array}%
	\right. ,\text{ \ \ }\alpha _{k}^{1}\equiv \alpha _{k}\equiv \left[ 
	\begin{array}{c}
	0 \\ 
	1 \\ 
	0 \\ 
	\vdots \\ 
	0%
	\end{array}%
	\right. ,\text{ \ \ \ }\alpha _{k}^{2}\equiv \alpha _{k}.\alpha _{k}\equiv %
	\left[ 
	\begin{array}{c}
	0 \\ 
	0 \\ 
	1 \\ 
	\vdots \\ 
	0%
	\end{array}%
	\right. ,\text{ \ \ \ \ \ }\alpha _{k}^{k}\equiv \alpha _{k}^{k-1}.\alpha
	_{k}\equiv \left[ 
	\begin{array}{c}
	0 \\ 
	0 \\ 
	0 \\ 
	\vdots \\ 
	1%
	\end{array}%
	\right.
	\end{equation*}%
	\begin{equation}
	a=\left[ 
	\begin{array}{c}
	a_{0} \\ 
	a_{1} \\ 
	\vdots \\ 
	a_{k}%
	\end{array}%
	\right. =a_{0}\left[ 
	\begin{array}{c}
	1 \\ 
	0 \\ 
	0 \\ 
	\vdots \\ 
	0%
	\end{array}%
	\right. +a_{1}\left[ 
	\begin{array}{c}
	0 \\ 
	1 \\ 
	0 \\ 
	\vdots \\ 
	0%
	\end{array}%
	\right. +\cdots a_{k}\left[ 
	\begin{array}{c}
	0 \\ 
	0 \\ 
	0 \\ 
	\vdots \\ 
	1%
	\end{array}%
	\right. =a_{0}\alpha _{k}^{0}+a_{1}\alpha _{k}^{1}+a_{2}\alpha
	_{k}^{2}+\cdots +a_{k}\alpha _{k}^{k}  \label{polyseriesform}
	\end{equation}
\end{definition}

\begin{algorithm}
	\text{ \ \ }\newline
	If%
	\begin{equation*}
	a=\left[ 
	\begin{array}{c}
	a_{0} \\ 
	a_{1} \\ 
	\vdots \\ 
	a_{k}%
	\end{array}%
	\right . \qquad \text{and}\ \ \ b=\left[ 
	\begin{array}{c}
	b_{0} \\ 
	b_{1} \\ 
	\vdots \\ 
	\ b_{k}%
	\end{array}%
	\right.,
	\end{equation*}
	then%
	\begin{equation*}
	a+b=\left[ 
	\begin{array}{c}
	a_{0}+b_{0} \\ 
	a_{1}+b_{1} \\ 
	\vdots \\ 
	a_{k}+\ b_{k}%
	\end{array}%
	\right . \ ;\ \ \ \ \ \ \ \lambda a=\left[ 
	\begin{array}{c}
	\lambda a_{0} \\ 
	\lambda a_{1} \\ 
	\vdots \\ 
	\lambda a_{k}%
	\end{array}%
	\right.
	\end{equation*}%
	\begin{equation}
	a\, . \, b=\left[ 
	\begin{array}{l}
	a_{0}b_{0} \\ 
	a_{1}b_{0}+a_{1}b_{0} \\ 
	a_{0}b_{2}+a_{1}b_{1}+a_{2}b_{0} \\ 
	\vdots \\ 
	a_{0}b_{k}+a_{1}b_{k-1}+\cdots +a_{k-1}b_{1}+a_{k}b_{0}%
	\end{array}%
	\right.  \label{polyseriesproduct}
	\end{equation}
\end{algorithm}

\begin{theorem}
	\text{ \ \ }\newline
	The identity \eqref{productpolyseriesidentity} holds for any truncated
	polyseries $a=\sum_{i=0}^{k}a_{i}\alpha _{k}^{i},$ and $b=%
	\sum_{i=0}^{k}b_{i}\alpha _{k}^{i}$ expressed in term of the $k-$basis $%
	\left[ \alpha _{k}^{0}\text{ \ }\alpha _{k}^{1}\text{ \ }\alpha _{k}^{2}%
	\text{ }\cdots \text{\ \ }\alpha _{k}^{k}\right] $ 
	\begin{equation}
	\left( \sum_{i=0}^{k}a_{i}\alpha _{k}^{i}\right)\, . \, \left(
	\sum_{i=0}^{k}b_{i}\alpha _{k}^{i}\right) =\sum_{i=0}^{k}\left(
	\sum_{j=0}^{i}a_{i}b_{j-i}\right) \alpha _{k}^{i}
	\label{productpolyseriesidentity}
	\end{equation}
\end{theorem}

\begin{proof}
	\text{ \ \ }\newline
	From \eqref{polyseriesform} we have:%
	\begin{eqnarray*}
		a &=&a_{0}\alpha _{k}^{0}+a_{1}\alpha _{k}^{1}+a_{2}\alpha _{k}^{2}+\cdots
		+a_{k}\alpha _{k}^{k}=\sum_{i=0}^{k}a_{i}\alpha _{k}^{i} \\
		&&\text{and} \\
		b &=&b_{0}\alpha _{k}^{0}+b_{1}\alpha _{k}^{1}+b_{2}\alpha _{k}^{2}+\cdots
		+b_{k}\alpha _{k}^{k}=\sum_{i=0}^{k}b_{i}\alpha _{k}^{i}
	\end{eqnarray*}
	Using \eqref{polyseriesproduct} we get:%
	\begin{equation*}
	\begin{array}{rl}
	a\, . \, b=\left( \sum_{i=0}^{k}a_{i}\alpha _{k}^{i}\right) \, . \,\left(
	\sum_{i=0}^{k}b_{i}\alpha _{k}^{i}\right) = & \left( a_{0}b_{0}\right)
	\alpha _{k}^{0}+ \\ 
	& +\left( a_{1}b_{0}+a_{1}b_{0}\right) \alpha _{k}^{1}+ \\ 
	& +\left( a_{0}b_{2}+a_{1}b_{1}+a_{2}b_{0}\right) \alpha _{k}^{2}+ \\ 
	& \vdots \\ 
	& +\left( a_{0}b_{k}+a_{1}b_{k-1}+\cdots +a_{k-1}b_{1}+a_{k}b_{0}\right)
	\alpha _{k}^{k} \\ 
	= & \sum_{i=0}^{k}\left( \sum_{j=0}^{i}a_{i}b_{j-i}\right) \alpha _{k}^{i}%
	\end{array}%
	\end{equation*}
\end{proof}

\section{Modular Resolution}

\begin{notation}
	\text{ \  \  }%
	
	\begin{description}
		\item[Integers and residues] If $m,$ $x\in \mathbb{Z}$, $y\in \mathbb{Z}$, we define:
		
		$x\equiv y\pmod m\ :\Longleftrightarrow \ m\mid
			(x-y):\Longleftrightarrow \exists k\in \mathbb{Z},$ $x-y=m\times
			k\Longleftrightarrow \left[ x\right] _{m}=\left[ y\right] _{m}$
			
			 $\left[ x\right] _{m}:=\left\{ t\in \mathbb{Z},\text{ }t\equiv x\pmod m%
			\right\} $
			
			 $\mathbb{Z}/m\mathbb{Z}:=\left\{ \left[ u\right] _{m},\text{ }u\in 
			\mathbb{Z}\right\} $

			\item[Modular Sets] 
		
		For any prime $p$, integers $a$,\,$b$,\,and sets $A$,\,$B$,\,(eventually finite as a fact)  we define:
		\begin{eqnarray*}
			\left[
			a,b\right]_{\mathbb{Z}} &:=&\left\{ c\in \mathbb{Z},\text{ }a\leq
			c\leq b\right\} \\
			\left( A \Subset B\text{ }\right) &:\Leftrightarrow & \bigg(\forall x\in \mathbb{Z},\,\text{ }
			\left(	x\in A\Big) \rightarrow \Big(\exists x^{\prime }\in B,\, x^{\prime }\equiv x\text{\ }\left( \mathrm{%
				mod}\text{ }p\right) \right)\bigg) \\
			\left( A \equiv B %
			\text{ }\right) &:\Leftrightarrow &\Big(\left( A \Subset B \right)\wedge \left( B \Subset A  \right)\Big)
		\end{eqnarray*}
	
		\item[Modular Fields] if $p$ is a prime number such that $p\geq 3,$ then:%
		\begin{eqnarray*}
			\mathbb{Z}/p\mathbb{Z} &\equiv &\mathbb{F}_{p}\equiv \left[ 0,\,p-1\right] _{%
				\mathbb{Z}} \\
			\mathbb{F}_{p}^{+} &\equiv &\left[ 0,\,\frac{p-1}{2}\right] _{\mathbb{Z}} \\
			\mathbb{F}_{p}^{-} &\equiv &\left[ \frac{p+1}{2},\,p-1\right] _{\mathbb{Z}%
			}\equiv \left[ \frac{1-p}{2},\,-1\right] _{\mathbb{Z}} \\
			0_{p} &\equiv &0_{\mathbb{F}_{p}} \\
			\left\vert m\right\vert _{p} &\equiv &\left\vert m\right\vert _{\mathbb{F}%
				_{p}}\equiv \left\{ 
			\begin{array}{cccccc}
				m &  &  & \text{if} &  & m\in \mathbb{F}_{p}^{+} \\ 
				&  &  &  &  &  \\ 
				p+m &  &  & \text{if} &  & m\in \mathbb{F}_{p}^{-}%
			\end{array}%
			\right. ,\text{ }m\in \mathbb{Z}\text{\ }
		\end{eqnarray*}%
		\item[Catalan Numbers] 
		\begin{equation*}
		C_{n}=\frac{1}{n+1}\,.\,\binom{2n}{n},\quad n=0,1,2,\ldots
		\end{equation*}%
		\begin{equation*}
		\begin{array}{l}
		C_{0}=1,\quad C_{1}=1,\quad C_{2}=2,\quad C_{3}=5,\quad C_{4}=14, \\ 
		C_{5}=42,\quad C_{6}=132,\quad C_{7}=429,\quad C_{8}=1430,\quad C_{9}=4862%
		\end{array}%
		\end{equation*}%
		\begin{equation}
		C_{0}=1,\quad C_{n+1}=\sum_{k=0}^{n}C_{k}\,.\,C_{n-k}=\frac{2(2n+1)}{n+2}%
		\,.\,C_{n},\quad (n\geq 0)  \label{CatalanFormula}
		\end{equation}%
		\begin{equation*}
		C=\frac{1-\left( 1-4u\right) ^{\frac{1}{2}}}{2u}=\sum_{n\geq
			0}C_{n}\,.\,u^{n}=1+u\,.\,C^{2},\quad u\in \left] 0,\frac{1}{4}\right[
		\end{equation*}

		\item[Polyseries] Set $u$ a polynumber, $\alpha _{1}$ and $\alpha _{2}$ two
		polyseries such that there exist polynumber sequences $\left( a_{k}\right)
		_{k\geq m},$ $\left( b_{k}\right) _{k\geq m},$ $m\in \mathbb{Z}$ verifying: 
		\begin{equation*}
		\alpha _{1}=\sum\limits_{k=m}^{+}a_{k}\cdot u^{k}\text{ \ and \ \ }\alpha
		_{2}=\sum\limits_{k=m}^{+}b_{k}\cdot u^{k}.
		\end{equation*}%
		We say that $\alpha _{1}$ is congruent to $\alpha _{2}$ modular $u^{n}$, $%
		n\in \mathbb{Z}$ if and only if there exists a polynumber sequence $c_{k},$
		such that: 
		\begin{equation*}
		\beta =\alpha _{1}-\alpha _{2}=\text{\ }\sum\limits_{k=n}^{+}c_{k}u^{k}.
		\end{equation*}%
		We write $\alpha _{1}\equiv \alpha _{2}\ \ \left( \mathrm{mod}\text{ }%
		u^{n}\right) $ and $\beta =0\left( u^{n}\right) $

	\end{description}
\end{notation}

\begin{proposition}
	\text{ \ \ }\newline
	For any positive constant integer $n\geq 3$ and any truncated polyserie $%
	A_{k}=\left[ \alpha _{0},\alpha _{1},\ldots \alpha _{k}\right\rangle $ , $%
	k\leq n,$ $\alpha _{k}\in \mathbb{F}_{p},$ if $u$ is a polynumber then we
	have: 
	\begin{equation}
	\left( \sum\limits_{k=0}^{n}\alpha _{k}\cdot u^{k}\right)
	^{2}=\sum\limits_{k=0}^{n}\left( \sum_{\ell =0}^{k}\alpha _{\ell }\,\cdot
	\alpha _{k-\ell }\right) u^{k}+\sum\limits_{k=n+1}^{2n}\left( \sum_{\ell
		=k-n}^{n}\alpha _{\ell }\,\cdot \alpha _{k-\ell }\right) \cdot u^{k}
	\label{squaresum}
	\end{equation}
\end{proposition}

\begin{proof}
	Set%
	\begin{equation*}
	S_{n}=\sum\limits_{k=0}^{n}\alpha _{k}u^{k}.
	\end{equation*}%
	We have:%
	\begin{equation}
	S_{n}^{2}=\sum_{k=0}^{2n}(\sum\limits_{\substack{ 0\leq i,j\leq n  \\ i+j=k}}%
	\alpha _{i}\alpha _{j})u^{k},  \label{squaresum1}
	\end{equation}
	
	with the equivalence of indices: 
	\begin{eqnarray*}
		\left( \ell =i\right) \wedge \left( i+j=k\right) &\Leftrightarrow &\left(
		\ell =i\right) \wedge \left( j=k-\ell \right) \\
		\left( 0\leq i\leq n\right) &\Leftrightarrow &\left( 0\leq \ell \leq n\right)
		\\
		\left( 0\leq j\leq n\right) &\Leftrightarrow &\left( 0\leq k-\ell \leq
		n\right) \Leftrightarrow \left( k-n\leq \ell \leq k\right) \\
		\left( 0\leq i\leq n\right) \wedge \left( 0\leq j\leq n\right)
		&\Leftrightarrow &\left( 0\leq \ell \leq n\right) \wedge \left( k-n\leq \ell
		\leq n\right) \\
		&\Leftrightarrow &\left( \max (0,\,k-n)\leq \ell \leq \min (k,\,n)\right) .
	\end{eqnarray*}
	
	Hence:%
	\begin{equation}
	\sum\limits_{\substack{ 0\leq i,j\leq n  \\ i+j=k}}\alpha _{i}\alpha
	_{j}=\sum_{\ell =\max (0,\,k-n)}^{\min (k\,,n)}\alpha _{\ell }\alpha
	_{k-\ell }.  \label{squaresum2}
	\end{equation}
	
	Using \eqref{squaresum2} we get:%
	\begin{equation*}
	\sum_{k=0}^{n}\left( \sum_{\ell =\max (0,\,k-n)}^{\min (k,\,n)}\alpha _{\ell
	}\alpha _{n-\ell }\right) u^{k}=\sum\limits_{k=0}^{n}\left( \sum_{\ell
		=0}^{k}\alpha _{\ell }\,\alpha _{k-\ell }\right) u^{k},
	\end{equation*}
	
	and%
	\begin{equation*}
	\sum_{k=n+1}^{2n}(\sum_{\ell =\max (0,\,k-n)}^{\min (k,\,n)}\alpha _{\ell
	}\alpha _{n-\ell })u^{k}=\sum\limits_{k=n+1}^{2n}(\sum_{\ell =k-n}^{n}\alpha
	_{\ell }\,\alpha _{k-\ell })u^{k}.
	\end{equation*}
	
	Substituting in \eqref{squaresum1} we get:%
	\begin{eqnarray*}
		S_{n}^{2} &=&\left( \sum\limits_{k=0}^{n}\alpha _{k}u^{k}\right)
		^{2}=\sum_{k=0}^{2n}(\sum\limits_{\substack{ 0\leq i,j\leq n  \\ i+j=k}}%
		\alpha _{i}\alpha _{j})u^{k}=\sum_{k=0}^{2n}(\sum_{\ell =\max
			(0,\,k-n)}^{\min (k,\,n)}\alpha _{\ell }\alpha _{n-\ell })u^{k} \\
		&=&\sum\limits_{k=0}^{n}(\sum_{\ell =0}^{k}\alpha _{\ell }\,\alpha _{k-\ell
		})u^{k}+\sum\limits_{k=n+1}^{2n}(\sum_{\ell =k-n}^{n}\alpha _{\ell }\,\alpha
		_{k-\ell })u^{k}
	\end{eqnarray*}
\end{proof}

\begin{remark}
	If $t\in \mathbb{F}_{p}^{+}$ and $u=t\ast e_{0}=t\ast \alpha ^{0}=\left[ 
	\begin{array}{cccc}
	t & 0 & 0 & \cdots%
	\end{array}%
	\right\rangle $ then the identity \eqref{squaresum} holds in $\mathbb{F}%
	_{p}^{+}.$
\end{remark}

\begin{proposition}
	\text{ \  \  }  \newline%
	{if }$n\in \mathbb{Z},\ n\geq 1,$ $u\in \mathbb{F}_{p}^{+},\ |u|\geq 2,$ $\
	t\in \mathbb{F}_{p}^{+},$ and defining $D_{u}$ as:\label{expansion}%
	\begin{equation*}
	D_{u}=\{0,1,\dots ,|u|-1\},
	\end{equation*}%
	then $\bigl|\{0,\dots ,|u|^{n}-1\}\bigr|\xrightarrow{\ \sim\ }\mathbb{Z}%
	/u^{n}\mathbb{Z},$ such that:%
	\begin{equation*}
	\exists !\ (a_{k})_{k=0}^{n-1}\in D_{u}^{\,n}:\quad t\equiv
	\sum_{k=0}^{n-1}a_{k}u^{k}\pmod{u^n}
	\end{equation*}
\end{proposition}

\begin{proof}	
	Set $r\equiv t\pmod{u^n},\ 0\leq r<|u|^{n},$ then:%
	\begin{eqnarray*}
		\exists q_{1},a_{0} &:&\ r=q_{1}u+a_{0},\ 0\leq a_{0}<|u| \\
		\exists q_{2},a_{1} &:&\ q_{1}=q_{2}u+a_{1},\ 0\leq a_{1}<|u| \\
		&&\vdots \\
		q_{n-1} &=&a_{n-1},\ 0\leq a_{n-1}<|u|
	\end{eqnarray*}%
	\begin{equation*}
	r=\sum_{k=0}^{n-1}a_{k}u^{k}\Rightarrow \bigl(\sum_{k=0}^{n-1}a_{k}u^{k}%
	\equiv \sum_{k=0}^{n-1}b_{k}u^{k}\pmod{u^n}\bigr)\Longrightarrow \
	\sum_{k=0}^{n-1}(a_{k}-b_{k})u^{k}\equiv 0\pmod{u^n}
	\end{equation*}
	
	If $m=\min \{k:a_{k}\neq b_{k}\},$ then:%
	\begin{equation*}
	u^{m}\Bigl(\sum_{j=0}^{n-1-m}(a_{m+j}-b_{m+j})u^{j}\Bigr)\equiv 0\pmod{u^n}
	\end{equation*}%
	\begin{eqnarray*}
		&\Rightarrow &\sum_{j=0}^{n-1-m}(a_{m+j}-b_{m+j})u^{j}\equiv 0%
		\pmod{u^{\,n-m}} \\
		&\Rightarrow &\Bigl|\sum_{j=0}^{n-1-m}(a_{m+j}-b_{m+j})u^{j}\Bigr|\leq
		(|u|-1)\sum_{j=0}^{n-1-m}|u|^{j}=|u|^{\,n-m}-1 \\
		&\Rightarrow &\sum_{j=0}^{n-1-m}(a_{m+j}-b_{m+j})u^{j}=0\Rightarrow \
		a_{k}=b_{k}
	\end{eqnarray*}
	
	We deduce that:%
	\begin{equation*}
	|u|\geq 2,\ n\geq 1\ \Longrightarrow \ \bigl\{0,\dots ,|u|^{n}-1\bigr\}%
	\xrightarrow{\ \sim\ }\mathbb{Z}/u^{n}\mathbb{Z}
	\end{equation*}
\end{proof}

\begin{lemma}
	\text{ \ \ }\newline
	\label{Catlan_Lem}{Set $p$}${>2}$ a prime{$,\ x\in \mathbb{F}_{p}^{+},$ $%
		n\in \mathbb{N},$ $n>2,$} then the equivalence \eqref{QuatEquaEqui} holds
	for any $t\in \mathbb{F}_{p}^{+}.$ 
	\begin{equation}
	\left( x^{2}-x+t\equiv 0\ \ \left( \mathrm{mod}\text{ }t^{n}\right) \right)
	\Leftrightarrow \left( x\equiv \sum\limits_{k=1}^{n-1}C_{k-1}\cdot \,t^{k}\
	\ \left( \mathrm{mod}\text{ }t^{n}\right) \right)  \label{QuatEquaEqui}
	\end{equation}
\end{lemma}

\begin{proof}
	\text{ \ \ }
	
	\begin{itemize}
		\item \textbf{\underline{The existence of The solution:}}\newline
		From proposition (\ref{expansion}), for any {$n\in \mathbb{N}_{{n>2}},$ and} 
		$x\in \mathbb{F}_{p}^{+},$ we have: 
		\begin{equation*}
		\exists \left\{ a_{k}\right\} _{k}\subset \mathbb{F}_{p}^{+}:x\equiv
		\sum\limits_{k=0}^{n-1}a_{k}\,t^{k}\text{ }\left( \mathrm{mod}\text{ }%
		t^{n}\right) ,
		\end{equation*}
		
		such that:%
		\begin{equation*}
		\left( x^{2}-x+t\equiv 0\ \ \left( \mathrm{mod}\text{ }t^{n}\right) \right)
		\Leftrightarrow \left( \sum\limits_{k=0}^{n-1}a_{k}t^{k}\right)
		^{2}-\sum\limits_{k=0}^{n-1}a_{k}t^{k}+t\ \equiv 0\ \ \left( \mathrm{mod}%
		\text{ }t^{n}\right) .
		\end{equation*}
		
		We have:%
		\begin{eqnarray*}
			\left( \sum\limits_{k=0}^{n-1}a_{k}\,t^{k}\right)
			^{2}-\sum\limits_{k=0}^{n-1}a_{k}\,t^{k}+t\ &=&\sum_{k=0}^{2\left(
				n-1\right) }(\sum\limits_{\substack{ 0\leq i,j\leq n  \\ i+j=k}}%
			a_{i}a_{j})t^{k}-\sum\limits_{k=0}^{n-1}a_{k}\,t^{k}+t \\
			&=&a_{0}-a_{0}^{2}+t+\sum\limits_{k=1}^{n-1}(a_{k}-\sum%
			\limits_{i+j=k-1}a_{i}a_{j})\,t^{k}.
		\end{eqnarray*}
		
		If $x^{2}-x+t\equiv 0\ \ \left( \mathrm{mod}\text{ }t^{n}\right) ,$ then: 
		\begin{eqnarray*}
			\left\{ 
			\begin{array}{c}
				a_{0}-a_{0}^{2}=0_{p} \\ 
				a_{1}-2a_{0}a_{1}+1=0_{p} \\ 
				a_{k}-\sum\limits_{i+j=k-1}a_{i}a_{j}=0_{p},\text{ \ \ }k>1%
			\end{array}%
			\right. &\Leftrightarrow &\left\{ 
			\begin{array}{l}
				a_{0}=a_{1}=1 \\ 
				a_{k}-\sum\limits_{i+j=k-1}a_{i}a_{j}=0_{p},\text{ \ \ }k>1%
			\end{array}%
			\right. \\
			&\Leftrightarrow &a_{k}=C_{k-1},\text{ }k>1,\,k\in 
			\mathbb{N}
			.
		\end{eqnarray*}
		
		We deduce that:%
		\begin{equation*}
		\left( x^{2}-x+t\equiv 0\ \ \left( \mathrm{mod}\text{ }t^{n}\right) \right)
		\Rightarrow x\equiv \sum\limits_{k=1}^{n-1}C_{k-1}\cdot \,t^{k}\ \ \left( 
		\mathrm{mod}\text{ }t^{n}\right) .
		\end{equation*}
		
		Suppose now that: 
		\begin{equation*}
		\ x\equiv C_{0}t+C_{1}t^{2}+C_{2}t^{3}+\cdots +C_{n-2}t^{n-1}\left( \mathrm{%
			mod}\text{ }t^{n}\right) =\sum\limits_{k=1}^{n-1}C_{k-1}\,t^{k}\left( 
		\mathrm{mod}\text{ }t^{n}\right) ,
		\end{equation*}%
		then: 
		\begin{equation}
		x=\sum\limits_{k=1}^{n-1}C_{k-1}\,t^{k}+\sum\limits_{k\geq
			n}a_{k}t^{k}=\alpha _{n}+\beta _{n}.  \label{CompleteForm}
		\end{equation}%
		So we get: 
		\begin{equation*}
		t^{2}=\left( \alpha _{n}+\beta _{n}\right) ^{2}=\alpha _{n}^{2}+\beta
		_{n}\left( 2\alpha _{n}+\beta _{n}\right) \equiv \alpha _{n}^{2}\text{ \ }(%
		\mathrm{mod}\ t^{n}).
		\end{equation*}%
		Since $\beta _{n}=\sum\limits_{k\geq n}a_{k}t^{k}\equiv 0$ $(\mathrm{mod}\
		u^{n}),$ we deduce that:%
		\begin{equation}
		x^{2}\equiv \alpha _{n}^{2}\ \ \left( \mathrm{mod}\text{ }t^{n}\right) .
		\label{quadeq1}
		\end{equation}%
		On the other hand, by the Catalan Formula \eqref{CatalanFormula} \ and
		applying the identity \eqref{squaresum} we get:%
		\begin{equation*}
		\alpha _{n}=\sum\limits_{k=0}^{\
			n}C_{k}\,t^{k+1}=\sum\limits_{k=1}^{n-1}C_{k-1}t^{k},
		\end{equation*}
		
		such that:%
		\begin{eqnarray*}
			\alpha _{n}^{2} &=&\left( \sum\limits_{k=1}^{n-1}C_{k-1}\,t^{k}\right)
			^{2}=t^{2}\left( \sum\limits_{i=0}^{n-1}C_{i}\,t^{i}\right) \left(
			\sum\limits_{j=0}^{n-1}C_{j}\,t^{j}\right) \\
			&=&u^{2}\left( \sum\limits_{k=0}^{n-1}\left( \sum_{\ell =0}^{k}C_{\ell
			}\,C_{k-\ell }\right) t^{k}+\sum\limits_{k=n}^{2\left( n-1\right) }\left(
			\sum_{\ell =k-n+1}^{n-1}C_{\ell }\,C_{k-\ell }\right) t^{k}\right)
		\end{eqnarray*}
		
		We deduce that:%
		\begin{equation*}
		\alpha _{n}^{2}\equiv t^{2}\sum\limits_{k=0}^{n-1}C_{k+1}\,t^{k}\text{ }(%
		\mathrm{mod}\ t^{n}).
		\end{equation*}
		
		We have%
		\begin{equation*}
		t^{2}\sum\limits_{k=0}^{n-1}C_{k+1}\,t^{k}=\sum\limits_{k=0}^{n-1}C_{k+1}%
		\,t^{k+2}=-t+\sum\limits_{k=0}^{\ n}C_{k}\,t^{k+1}\equiv
		-t+\sum\limits_{k=1}^{n-1}C_{k-1}\,t^{k}\text{ }(\mathrm{mod}\ t^{n}).
		\end{equation*}%
		Then:%
		\begin{equation}
		\alpha _{n}^{2}\equiv -t+x\ \ \left( \mathrm{mod}\text{ }t^{n}\right) .
		\label{quadeq2}
		\end{equation}
		
		From\ \eqref{quadeq1} and \eqref{quadeq2}\ we\ get $\ x^{2}\equiv -t+x$ $\
		\left( \mathrm{mod}\text{ }t^{n}\right) .$
		
		Hence:%
		\begin{equation}
		x^{2}-x+t\equiv 0\ \ \left( \mathrm{mod}\text{ }t^{n}\right) .
		\label{ModEqua1}
		\end{equation}
		
		\item \textbf{\underline{The Uniqueness of the solution:}}
		
		\text{\ }
		 \newline
		 \vspace{20pt}
		Set $y\equiv x^{2}-x+t$ and $\ y^{\prime }\equiv x^{2}-x+t,$ then: 
		\begin{eqnarray*}
			\left\vert y-y^{\prime }\right\vert _{p} &\equiv &\left\vert x^{2}-x+t-x^{\prime }{%
				^{2}}+x^{\prime }-t\right\vert _{p} \\
			&\equiv &\left\vert \left( x-x^{\prime }\right) \left( x+x^{\prime }-1\right)
			\right\vert _{p}.
		\end{eqnarray*}%
		Since $x\in \mathbb{F}_{p}^{+}$ and $x^{\prime }\in \mathbb{F}_{p}^{+},$
		then: 
		\begin{equation*}
		\begin{tabular}{lll}
		$x,$ $x^{\prime }\in \left[ 1,\displaystyle\frac{p-1}{2}\right]_{\mathbb{Z}} $ & $%
		\Rightarrow $ & $x+x^{\prime }\in \left[ 2,p-1\right]_{\mathbb{Z}} $ \\ 
		& $\Rightarrow $ & $\left( x+x^{\prime }-1\right) \in \left[ 1,p-2\right]_{\mathbb{Z}} $
		\\ 
		& $\Rightarrow $ & $\left( x+x^{\prime }-1\right) \in \mathbb{F}_{p}^{\times
		}$ \\ 
		& $\Rightarrow $ & $\left( x+x^{\prime }-1\right) \not\equiv 0_{p}.$ \\ 
		& $\Rightarrow $ & $\left( x+x^{\prime }-1\right) \not\equiv 0$ $\ \left( 
		\mathrm{mod}\text{ }p\right) .$%
		\end{tabular}%
		\end{equation*}%
		We have the \ following equivalences: 
		\begin{equation*}
		\begin{tabular}{lll}
		$\left\vert y-y^{\prime }\right\vert _{p}\equiv 0_{p}$ & $\Leftrightarrow $ & $%
		\left\vert \left( x-x^{\prime }\right) \left( x+x^{\prime }-1\right)
		\right\vert _{p}\equiv 0_{p}$ \\ 
		& $\Leftrightarrow $ & $\left\vert x-x^{\prime }\right\vert _{p}\equiv 0_{p}$ \\ 
		& $\Leftrightarrow $ & $x-x^{\prime }\equiv 0_{p}$  \\ 
		& $\Leftrightarrow $ & $x$ \ $\equiv x^{\prime }$ $\left( \mathrm{mod}\text{ 
		}p\right) $%
		\end{tabular}%
		\end{equation*}%
		The uniqueness of the solution of the modular quadratic equation %
		\eqref{ModEqua1} holds.
	\end{itemize}
\end{proof}

\begin{remark}
	If we take $t\equiv \displaystyle\frac{1}{2^{k}}$ $\left( \mathrm{mod}\text{ 
	}p\right) ,\,(k>2),$ we can apply the Catalan Formula \eqref{CatalanFormula}%
	, since the dyadic values of $t$ assures that $t\in \left] 0,\displaystyle%
	\frac{1}{4}\right[ _{\mathbb{F}_{p}}.$
\end{remark}

\begin{example}
	Set $p=11,$ we get: 
	\begin{equation*}
	\ \forall q\in \mathbb{Z}_{11}^{\times }:\ \mathrm{ord}_{11}(q)\;\bigm|%
	\;\varphi (11)=11-1=10,
	\end{equation*}%
	such that:%
	\begin{equation*}
	\forall q\in \mathbb{Z}_{11}^{\times },\;\exists k\in \mathbb{N}_{+}:\;k=%
	\mathrm{ord}_{11}(q)\quad \text{and}\quad q^{k}\equiv 1\pmod{11}.
	\end{equation*}
	
	$\mathbb{F}_{11}^{+}\equiv\left[ 1,\displaystyle\frac{11-1}{2}\right]_{\mathbb{Z}} \equiv\left[ 1,5%
	\right]_{\mathbb{Z}} ,$ and $q=2$ $\Rightarrow \ q^{5}\equiv -1\pmod{11}$ $\Rightarrow \
	\ q^{10}\equiv +1\pmod{11}.$ \text{ \ \ } \newline
	\text{ \ \ } \newline
	We have: \ 
	\begin{equation*}
	 ]0,\frac{1}{4}[_{\mathbb{F}_{11}}\text{ }\equiv \text{ }%
	\{2^{-10},2^{-9},2^{-8},2^{-7},2^{-6},2^{-5},2^{-4},2^{-3}\}\text{ }%
	_{\mathbb{F}_{11}}\equiv\{1,2,4,8,5,10,9,7\}\text{ }_{\mathbb{F}_{11}}.
	\end{equation*}
	
	Taking for example $t\equiv\displaystyle\frac{1}{2^{9}}$ and $n=4>2,$ we get:%
	\begin{equation*}
	t\equiv2^{-9}\equiv 2\text{ }\left( \mathrm{mod}\text{ }11\right) \Rightarrow
	t^{n}\equiv 2^{4}\equiv5\text{ }%
	\left( \mathrm{mod}\text{ }11\right) .
	\end{equation*}%
	We can find the solution $x$ in $\left[ 1,5\right] $ :%
	\begin{equation*}
	\begin{tabular}{lll}
	$x^{2}-x+t\equiv 0\ \ \left( \mathrm{mod}\text{ }t^{n}\right) $ & $%
	\Longleftrightarrow $ & $x^{2}-x+2\equiv 0\ \text{ }\left( \mathrm{mod}\text{
	}5\right) $ \\ 
	& $\Longleftrightarrow $ & $x\equiv1\text{ }\left( \mathrm{mod}\text{ }5\right) .$%
	\end{tabular}%
	\end{equation*}
	
	We check that:%
	\begin{equation*}
	\sum\limits_{k=1}^{n-1}C_{k-1}\,\cdot t^{k}\equiv   C_{0}+C_{1}\cdot t+C_{2}\cdot
	t^{2}+C_{3}\cdot t^{3}\equiv \left( 1+2+2\cdot 2^{2}+5\cdot 2^{3}\right) 
	\text{ }\left( \mathrm{mod}\text{ }5\right) \equiv 1 \equiv x\text{ }\left( \mathrm{mod}\text{ }t^{n}\right)
	\end{equation*}
\end{example}
\begin{theorem}
	\label{Catlan_Thm}{Set $p$}${>2}$ a prime{\ and $\ x,a,t\in \mathbb{F}_{p}$
		such that $n\in \mathbb{N}$} {for $n>2,$} then the equivalence %
	\eqref{QuatEquaParEqui} holds 
	\begin{equation}
	\left( x^{2}+a\cdot x+t\equiv \text{ }0\text{ }\left( \mathrm{mod}\text{ }%
	t^{n}\right) \right) \Leftrightarrow \left\{ 
	\begin{array}{l}
	x=x_{1}\equiv \sum\limits_{k=1}^{n-1}C_{k-1}\cdot \displaystyle\frac{1}{%
		\left( -a\right) ^{2k-1}}\cdot t^{k}\text{ \ }\left( \mathrm{mod}\text{ }%
	t^{n}\right) \\ 
	\text{ or} \\ 
	x=x_{2}=-a-x_{1}\equiv -a-\sum\limits_{k=1}^{n-1}C_{k-1}\displaystyle\cdot 
	\frac{1}{\left( -a\right) ^{2k-1}}\cdot t^{k}\text{ }\left( \mathrm{mod}%
	\text{ }t^{n}\right)%
	\end{array}%
	\right. \text{ }  \label{QuatEquaParEqui}
	\end{equation}
\end{theorem}

\begin{proof}
	To find the solutions of the modular quadratic equation \eqref{QuatEquaPar}: 
	\begin{equation}
	x^{2}+a\cdot x+t\equiv \text{ }0\text{ }\left( \mathrm{mod}\text{ }%
	t^{n}\right) ,  \label{QuatEquaPar}
	\end{equation}%
	we transform it to the usual form \eqref{ModEqua1} 
	\begin{equation*}
	\begin{tabular}{lll}
	$x^{2}+a\cdot x+t\equiv $ $0$ $\left( \mathrm{mod}\text{ }t^{n}\right) $ & $%
	\Leftrightarrow $ & $\displaystyle\left( \frac{x}{-a}\right) ^{2}-\left( 
	\frac{x}{-a}\right) +\frac{t}{a^{2}}\equiv \text{ }0\text{ }\left( \mathrm{%
		mod}\text{ }t^{n}\right) $ \\ 
	&  &  \\ 
	& $\Leftrightarrow $ & $y^{2}-y+t\equiv 0$ $\left( \mathrm{mod}\text{ }%
	t^{n}\right) .$%
	\end{tabular}%
	\end{equation*}
	
	\textbf{case 1: }$\ {x=x_{1}}\in \mathbb{F}_{p}^{+}$\newline
	
	From Lemma \eqref{Catlan_Lem} {\ }we get:%
	\begin{equation*}
	\frac{x_{1}}{-a}=y\equiv \alpha _{n}=\sum\limits_{k=1}^{n-1}C_{k-1}\,t^{k}\
	\ =\sum\limits_{k=1}^{n-1}C_{k-1}\,\left( \frac{t}{a^{2}}\right) ^{k}\ \text{
	}\left( \mathrm{mod}\text{ }t^{n}\right) \ .
	\end{equation*}%
	$\Rightarrow $ 
	\begin{equation*}
	x_{1}\equiv -a\sum\limits_{k=1}^{n-1}C_{k-1}\,\frac{t^{k}}{a^{2k}}\ \
	=\sum\limits_{k=1}^{n-1}C_{k-1}\,\frac{t^{k}}{\left( -a\right) ^{2k-1}}\
	\left( \mathrm{mod}\text{ }t^{n}\right) .\ 
	\end{equation*}%
	We deduce that\setlength{\abovedisplayskip}{15pt} 
	\begin{equation}
	x_{1}\equiv \sum\limits_{k=1}^{n-1}C_{k-1}\frac{1}{\left( -a\right) ^{2k-1}}%
	t^{k}\text{ \ }\left( \mathrm{mod}\text{ }t^{n}\right) .  \label{Sol1}
	\end{equation}
	
	\textbf{case 2: }{{$x=x_{2}$}$\in \mathbb{F}_{p}^{-}$}\newline
	
	We have: $\ $%
	\begin{eqnarray*}
		\left( -a-x_{2}\right) ^{2}+a\left( -a-x_{2}\right) +t &=&a^{2}+2a\cdot
		x_{2}+x_{2}^{2}-a^{2}-a\cdot x_{2}+t \\
		&=&x_{2}^{2}+a\cdot x_{2}+t
	\end{eqnarray*}
	
	Assuming that $x_{1}=-a-x_{2},$ and since {\ $x_{2}\in \mathbb{F}_{p}^{-},$
		then } $x_{1}\in \mathbb{F}_{p}^{+}.$
	
	By using \eqref{Sol1} we get: 
	\begin{equation*}
	-a-x_{2}=x_{1}\equiv \sum\limits_{k=1}^{n-1}C_{k-1}\frac{1}{\left( -a\right)
		^{2k-1}}t^{k}\text{ \ }\left( \mathrm{mod}\text{ }t^{n}\right) .
	\end{equation*}
	
	Hence:%
	\begin{equation*}
	x_{2}\equiv -a-\sum\limits_{k=1}^{n-1}C_{k-1}\frac{1}{\left( -a\right)
		^{2k-1}}t^{k}\text{ \ }\left( \mathrm{mod}\text{ }t^{n}\right) =-a-x_{1}%
	\text{ }\left( \mathrm{mod}\text{ }t^{n}\right) .
	\end{equation*}
\end{proof}

\begin{example}
	If $a=+1,$ and $n=10,$ then the two solutions of the modular equation 
	\begin{equation*}
	x^{2}+x+t\equiv 0\text{ }\left( \mathrm{mod}\text{ }t^{10}\right) ,
	\end{equation*}%
	are: 
	\begin{equation*}
	x_{1}\equiv
	-t-t^{2}-2t^{3}-5t^{4}-14t^{5}-42t^{6}-132t^{7}-429t^{8}-1430t^{9}\text{ }%
	\left( \mathrm{mod}\text{ }t^{10}\right) ,
	\end{equation*}%
	and its modular conjuagate
	\begin{equation*}
	x_{2}\equiv
	-1+t+t^{2}+2t^{3}+5t^{4}+14t^{5}+42t^{6}+132t^{7}+429t^{8}+1430t^{9}\text{ }%
	\left( \mathrm{mod}\text{ }t^{10}\right) .
	\end{equation*}
	
	Up to degree $10,$ we get 
	\begin{equation*}
	\ \left( x_{1}^{2}+x_{1}+t\right) \text{ }\left( \mathrm{mod}\text{ }%
	t^{10}\right) =\left( x_{2}^{2}+x_{2}+t\right) \text{ }\left( \mathrm{mod}%
	\text{ }t^{10}\right) =0\text{ }\left( \mathrm{mod}\text{ }t^{10}\right) ,
	\end{equation*}%
	such that: 
	\begin{equation*}
	\resizebox{1.1\textwidth}{!}{$ \beta _{n}\left( t\right) =t^{10}\left(
		2044900\,t^{8}+1226940\,t^{7}+561561\,t^{6}+233376\,t^{5}+93500\,t^{4}+37400%
		\,t^{3}+15470\,t^{2}+7072\,t+4862\right) $}.
	\end{equation*}%
	For more complete examples \textbf{See} \cite{wildberger2025hyper}
\end{example}
\pagestyle{plain}
\normalsize 
\section{Conclusion}
The equivalence \ \eqref{QuatEquaParEqui} in the Theorem \eqref{Catlan_Thm},
shows the importance of the algebraic notion of truncation of polyseries to
prove the open problems about the existence and the uniqueness of solutions
of some kinds of modular quadratic equations over finite fields.
 \vspace{1.5cm}
\appendix
\begin{center}
	{ {\Large {{\textbf{\underline{Appendices}}}}} }
\end{center}
\appendix
\makeatletter
\let\orig@addcontentsline 
\addcontentsline
\renewcommand{\addcontentsline}[3]{}   \makeatother
 \vspace{1.5cm}
 \section{Algebra}
 
 \begin{definition}
 	A unitary commutative ring $\left( \mathrm{R},+,\cdot \right) $ with
 	addition unit $0_{\mathrm{R}}$ and \ multiplication unit $1_{\mathrm{R}}$ is
 	the algebraic structure verifying the following:
 	
 	\begin{description}
 		\item[additive associativity] $\forall \,x,y,z\in \mathrm{R}%
 		,\quad\;(x+y)+z=x+(y+z)$
 		
 		\item[additivec commutativity] $\forall \,x,y\in \mathrm{R}:\;x+y=y+x$
 		
 		\item[additive identity element] $\exists \,0_{\mathrm{R}}\in \mathrm{R}%
 		,\quad\;\forall \,x\in \mathrm{R}:\;x+0_{\mathrm{R}}=x$
 		
 		\item[additive inverse element] $\forall \,x\in R:\;\exists \,y\in \mathrm{R}%
 		,\quad\;x+y=0_{\mathrm{R}}$
 		
 		\item[Associative Multiplication] $\forall \,x,y,z\in \mathrm{R}:\;(x\cdot
 		y)\cdot z=x\cdot (y\cdot z)$
 		
 		\item[multiplicative identity element] $\exists \,1_{\mathrm{R}}\in \mathrm{R%
 		},\quad\;\forall \,x\in \mathrm{R},\quad\;1_{\mathrm{R}}\cdot x=x$
 		
 		\item[multiplicative commutativity] $\forall \,x,y\in \mathrm{R}:\;x\cdot
 		y=y\cdot x$
 		
 		\item[multiplicative associativity] $\forall \,x,y,z\in \mathrm{R}:\;x\cdot
 		(y+z)=x\cdot y+x\cdot z$
 		
 		\item[Unit Group of multiplicative inverse elements] The unit group $\mathbf{%
 			U}(\mathrm{R})$ elements is defined by:%
 		\begin{equation*}
 		\mathbf{U}(\mathrm{R})=\{\,u\in \mathrm{R}:\exists \,v\in \mathrm{R}%
 		:\;u\cdot v=v\cdot u=1_{\mathrm{R}}\,\}\neq \varnothing
 		\end{equation*}
 			\item[Finite field and polynomial ring] If $q=p^{r},\ p\ $prime$,\ r\in 
 		\mathbb{Z}_{\geq 1},$ $\mathbb{F}_{q}\ $(field)$,\ |\mathbb{F}_{q}|=q,$ then:%
 		\begin{equation*}
 		\mathbb{F}_{q}[t]=\left\{ \sum_{i=0}^{d}a_{i}t^{i}:d\geq 0,\ a_{i}\in 
 		\mathbb{F}_{q}\right\}
 		\end{equation*}
 		
 		\item[Principal Ideal Domain property] If $\deg :\mathbb{F}_{q}[t]\setminus
 		\{0\}\rightarrow \mathbb{Z}_{\geq 0},$ such that:%
 		\begin{equation*}
 		\forall f,g\in \mathbb{F}_{q}[t],\ g\neq 0,\ \exists q,r\in \mathbb{F}%
 		_{q}[t]:\ f=qg+r,\ (r=0\ \vee \ \deg r<\deg g)
 		\end{equation*}
 		
 		then $\mathbb{F}_{q}[t]\ $is Euclidean.
 		
 		\item[Irreducibles and residue fields] If \ $f\in \mathbb{F}_{q}[t],\ f\ $is
 		irreducible,$\ $such that $\deg f=d,$ then:%
 		\begin{equation*}
 		\mathbb{F}_{q}[t]/(f)\cong \mathbb{F}_{q^{d}},\ |\mathbb{F}_{q}[t]/(f)|=q^{d}
 		\end{equation*}
 		
 		\item[Valuations] If $\ f\ $irreducible, $\ $such that $v_{f}:\mathbb{F}%
 		_{q}(t)^{\times }\rightarrow \mathbb{Z},$ then:%
 		\begin{eqnarray*}
 			v_{f}(g/h) &=&{ord}_{f}(g)-{ord}_{f}(h) \\
 			v(g/h) &=&\deg (h)-\deg (g)
 		\end{eqnarray*}
 		
 		\item[Basis representation] If $b\in \mathbb{Z},\ b>1,$ then:%
 		\begin{equation*}
 		\forall n\in \mathbb{Z}_{\geq 1}\ \exists !\ (r_{i})_{i=0}^{k},\ r_{i}\in
 		\{0,\dots ,b-1\},\ r_{k}\neq 0:\ n=\sum_{i=0}^{k}r_{i}b^{i}
 		\end{equation*}
 		
 		If $p\in \mathbb{Z},\ p\geq 2,$ $x\in \lbrack 0,1),$ then:%
 		\begin{equation*}
 		\exists (a_{n})_{n\geq 1},\ a_{n}\in \{0,\dots ,p-1\}:\ x=\sum_{n\geq
 			1}a_{n}p^{-n}
 		\end{equation*}
 		
 		such that $\ a_{n}=\lfloor pr_{n-1}\rfloor ,\ r_{0}=x,\
 		r_{n}=pr_{n-1}-a_{n},\ 0\leq r_{n}<1,$ with:%
 		\begin{equation*}
 		S_{n}=\sum_{i=1}^{n}a_{i}p^{-i}\Rightarrow x-S_{n}=p^{-n}r_{n}
 		\end{equation*}
 		
 		and 
 		\begin{equation*}
 		\left( \sum_{n=1}^{+}a_{n}p^{-n}=\sum_{n=1}^{+}b_{n}p^{-n}\right)
 		\Leftrightarrow \left( (a_{n})=(b_{n})\ \vee \ \exists k:\ 
 		\begin{cases}
 		a_{k}=b_{k}+1 \\ 
 		\forall n>k:\ a_{n}=0,\ b_{n}=p-1%
 		\end{cases}%
 		\right)
 		\end{equation*}
 	\end{description}
 \end{definition}
 
 \begin{example}
 	\text{ \ \ }
 	
 	\begin{itemize}
 		\item $\mathbf{U}(\mathbb{Z}_{n})=\mathbb{Z}_{n}^{\times }$ $\neq
 		\varnothing $
 		
 		\item $\mathbf{U}(\mathbb{Z}_{10})=\left\{ 1,3,7,9\right\} $ since $%
 		3\,.\,7=21=1+2\cdot 10$ and $9\cdot 9=81=1+8\cdot 10$
 		
 		\item $\mathbf{U}(\mathbb{Z}_{11})=\left\{ 1,2,3,4,5,6,7,8,9\right\} =%
 		\mathbb{Z}_{11}-\left\{ 0\right\} $
 	\end{itemize}
 \end{example}
 \begin{corollary}
 	For any ring $\mathrm{R}$ and $\ x\in \mathbf{U}(\mathrm{R})$ the following
 	equivalences hold%
 	\begin{equation*}
 	\bigl(x^{2}-s\cdot \,x+p\in \mathrm{R}\bigr)\Leftrightarrow \bigl(s\in 
 	\mathrm{R}\wedge p\in \mathrm{R}\bigr)
 	\end{equation*}
 	
 	and%
 	\begin{equation*}
 	\bigl(x^{2}-s\,\cdot x+p=0_{\mathrm{R}}\bigr)\;\Leftrightarrow \bigl(%
 	\;\exists \,\alpha ,\beta \in \mathrm{R}:\left( \beta =s\right) \wedge
 	\left( p=\alpha \cdot \left( \beta -\alpha \right) \right) \wedge \left(
 	x\in \left\{ \alpha ,\beta -\alpha \right\} \right) \bigr)
 	\end{equation*}
 \end{corollary}
 \begin{proof}
 	\text{ \ \ }\newline
 	1) To prove the first equivalence, we prove the implication: 
 	\begin{eqnarray*}
 		\left( \left( x^{2}-s\cdot \,x+p\in \mathrm{R}\right) \wedge \left( x\in 
 		\mathbf{U}(\mathrm{R})\right) \right) &\Rightarrow &(x^{2}-s\,\cdot
 		x+p)\,x^{-1}\;=\;x\;-\;s\;+\;p\cdot \,x^{-1}\in \mathrm{R} \\
 		&\Rightarrow &s\;=\;x\;+\;p\cdot \,x^{-1}\;-\;\bigl((x^{2}-s\,\cdot
 		x+p)\,\cdot x^{-1}\bigr)\in \mathrm{R} \\
 		&\Rightarrow &p\;=\;x\cdot \,s\;-\;x^{2}\;\in \;\mathrm{R}.
 	\end{eqnarray*}
 	
 	The inverse implication is evident by the closure under multiplication,
 	additon and subtraction.
 	
 	2) Assume $\left( x^{2}-s\cdot \,x+p\;=\;0\right) \wedge \left( x\in \mathbf{%
 		U}(\mathrm{R})\right) ,$ then:%
 	\begin{eqnarray*}
 		x^{2}-s\cdot \,x+p\; &=&0_{\mathrm{R}}\Rightarrow p\;=\;s\,\cdot x\;-\;x^{2}
 		\\
 		\left( x=\alpha \right) \wedge \left( \beta \;=\;s\right) &\Rightarrow
 		&p=\alpha \,\cdot \left( \beta -\alpha \right) =x\cdot \;\left( s-x\right)
 		=\;x\,\cdot s-x^{2}\;\; \\
 		&\Rightarrow &\exists \,\alpha ,\beta \in \mathrm{R}:\;(\beta =s)\wedge
 		(p=\alpha \cdot \,\left( \beta -\alpha \right) ). \\
 		\left( x=\beta -\alpha \right) \wedge \left( \beta \;=\;s\right)
 		&\Rightarrow &p=\alpha \,\cdot \left( \beta -\alpha \right) =\left(
 		s-x\right) \cdot \left( x\right) =\;x\,\cdot s-x^{2}\;\; \\
 		&\Rightarrow &\exists \,\alpha ,\beta \in \mathrm{R}:\;(\beta =s)\wedge
 		(p=\alpha \cdot \,\left( \beta -\alpha \right) ).
 	\end{eqnarray*}
 	
 	Hence:%
 	\begin{equation*}
 	\bigl(x^{2}-s\,\cdot x+p=0_{\mathrm{R}}\bigr)\;\Leftrightarrow \bigl(%
 	\;\exists \,\alpha ,\beta \in \mathrm{R}:\left( \beta =s\right) \wedge
 	\left( p=\alpha \cdot \left( \beta -\alpha \right) \right) \wedge \left(
 	x\in \left\{ \alpha ,\beta -\alpha \right\} \right) \bigr)
 	\end{equation*}
 \end{proof}
 
 \begin{example}
 	Solving $x^{2}-3x+2=0$
 	
 	\begin{itemize}
 		\item In $\mathbb{Z}_{10}$ we have $x\in \left\{ 1\right\} $ since $\beta
 		=s=3\in \mathbf{U}(\mathbb{Z}_{10})$ and $\beta -\alpha =2\not\in \mathbf{U}(%
 		\mathbb{Z}_{10})$
 		
 		\item In $\mathbb{Z}_{11}$ we have $x\in \left\{ 1,2\right\} $ since $\beta
 		=s=3\in \mathbf{U}(\mathbb{Z}_{11})$ and $\beta -\alpha =2\in \mathbf{U}(%
 		\mathbb{Z}_{11})$
 	\end{itemize}
 \end{example}
 
\begin{description}
	\item[\protect\underline{Definitions}]\text{\ \ } 
	
	\begin{itemize}
		\item 
		$R$ commutative ring, $1_{R}\neq 0$,\, $k\in \mathbb{N}_{>1}$
		\item $\mathcal{A}_{k}(R)=\{(a_{0},a_{1},\dots ,a_{k})\mid a_{i}\in R\}$.
		
		\item $(a+b)_{i}=a_{i}+b_{i}$, $0\leq i\leq k$.
		
		\item $(a\ast b)_{n}=\sum_{i=0}^{n}a_{i}b_{n-i}$, $0\leq n\leq k$.
		
		\item $(r\cdot a)_{i}=ra_{i}$, $r\in R$.
		
		\item $0=(0,0,\dots ,0)$, $1=(1_{R},0,\dots ,0)$.
	\end{itemize}
	\vspace{0.7em}
	\item[\protect\underline{Module Axioms}]\text{\ \ }\\ 
	$\because $ 
	\begin{itemize}
		\item  $(a,b,c\in \mathcal{A}%
		_{k}(R),r,s\in R)$
		\item $%
		((a+b)+c)_{i}=(a+b)_{i}+c_{i}=(a_{i}+b_{i})+c_{i}=a_{i}+(b_{i}+c_{i})=a_{i}+(b+c)_{i}=(a+(b+c))_{i}
		$
		
		\item $(a+b)_{i}=a_{i}+b_{i}=b_{i}+a_{i}=(b+a)_{i}$
		
		\item $(a+0)_{i}=a_{i}+0=a_{i}=0+a_{i}=(0+a)_{i}$
		
		\item $(-a)_{i}=-a_{i}$, $%
		(a+(-a))_{i}=a_{i}+(-a_{i})=0=(-a)_{i}+a_{i}=((-a)+a)_{i}$
		
		\item $(r\cdot (a+b))_{i}=r(a_{i}+b_{i})=ra_{i}+rb_{i}=(r\cdot
		a)_{i}+(r\cdot b)_{i}=(r\cdot a+r\cdot b)_{i}$
		
		\item $((r+s)\cdot a)_{i}=(r+s)a_{i}=ra_{i}+sa_{i}=(r\cdot a)_{i}+(s\cdot
		a)_{i}=(r\cdot a+s\cdot a)_{i}$
		
		\item $((rs)\cdot a)_{i}=(rs)a_{i}=r(sa_{i})=r\cdot (s\cdot a)_{i}=(r\cdot
		(s\cdot a))_{i}$
		
		\item $(1_{R}\cdot a)_{i}=1_{R}\cdot a_{i}=a_{i}$
	\end{itemize}
	\vspace{20pt}
	$\therefore $\, $(\mathcal{A}_{k}(R),+)
	$ abelian group.
	\vspace{0.7em}
	\item[\protect\underline{Associativity of $\ast $}] 
	\text{  \  }\\ %
	
	$\because $ $ a,b,c\in \mathcal{A}_{k}(R),\,0\leq n\leq k,\,$ $n\in 
	\mathbb{N}
	_{>1}$

	$\because $ \ $\ \left( \left( j=m-i\right) \wedge \left( \ell
	=n-i-j\right) \right) \Leftrightarrow \left( \left( m=i+j\right) \wedge
	\left( n=i+j+\ell \right) \right) $
	
	$\therefore $ $\ $%
	\begin{eqnarray*}
		((a\ast b)\ast c)_{n} &=&\sum_{m=0}^{n}(a\ast
		b)_{m}c_{n-m}=\sum_{m=0}^{n}\left( \sum_{i=0}^{m}a_{i}b_{m-i}\right)
		c_{n-m}=\sum_{m=0}^{n}\sum_{i=0}^{m}a_{i}b_{m-i}c_{n-m} \\
		&=&\sum_{i=0}^{n}\sum_{j=0}^{n-i}a_{i}b_{j}c_{n-i-j}=\sum_{\substack{ %
				i+j+\ell =n \\ i,j,\ell \geq 0}}a_{i}b_{j}c_{\ell }
	\end{eqnarray*}
	
	$\therefore $%
	\begin{eqnarray*}
		(a\ast (b\ast c))_{n} &=&\sum_{i=0}^{n}a_{i}(b\ast
		c)_{n-i}=\sum_{i=0}^{n}a_{i}\left( \sum_{j=0}^{n-i}b_{j}c_{n-i-j}\right)
		=\sum_{i=0}^{n}\sum_{j=0}^{n-i}a_{i}b_{j}c_{n-i-j} \\
		&=&\sum_{\substack{ i+j+\ell =n \\ i,j,\ell \geq 0}}a_{i}b_{j}c_{\ell
		}=((a\ast b)\ast c)_{n}
	\end{eqnarray*}
	
	$\because $%
	\begin{equation*}
	(a\ast (b\ast c))_{n}=((a\ast b)\ast c)_{n},\left( n\leq k\right) 
	\end{equation*}
	
	$\therefore $%
	\begin{equation*}
	a\ast (b\ast c))= (a\ast b)\ast c,\ ( a,b,c\in \mathcal{A}_{k}(R))
	\end{equation*}
	
	$\therefore $ $\mathbf{associative}\left( \mathcal{A}_{k}(R),\ast \right) $
	\vspace{0.7em}
	
	\item[\protect\underline{Left Distributivity}] \text{\ }
	\\
	$\because $%
	\begin{eqnarray*}
		(a\ast (b+c))_{n}
		&=&\sum_{i=0}^{n}a_{i}(b+c)_{n-i}=\sum_{i=0}^{n}a_{i}(b_{n-i}+c_{n-i})=%
		\sum_{i=0}^{n}(a_{i}b_{n-i}+a_{i}c_{n-i}) \\
		&=&\sum_{i=0}^{n}a_{i}b_{n-i}+\sum_{i=0}^{n}a_{i}c_{n-i}=(a\ast
		b)_{n}+(a\ast c)_{n}=((a\ast b)+(a\ast c))_{n}
	\end{eqnarray*}
	
	$\because $%
	\begin{equation*}
	(a\ast (b+c))_{n}=(a\ast (b+c))_{n},\left( n\leq k\right) 
	\end{equation*}

	$\therefore $%
	\begin{equation*}
	a\ast (b+c)=a\ast (b+c),\ (a,b,c\in \mathcal{A}_{k}(R)
	\end{equation*}
	
	$\therefore $ $\mathbf{Left\_Distributivity}\left( \mathcal{A}_{k}(R),\ast
	,.\right) $
	
	\vspace{0.7em}
	\item[\protect\underline{Right Distributivity}]\text{\ }
	\\
	$\because $%
	\begin{eqnarray*}
		((b+c)\ast a)_{n}
		&=&\sum_{i=0}^{n}(b+c)_{i}a_{n-i}=\sum_{i=0}^{n}(b_{i}+c_{i})a_{n-i}=%
		\sum_{i=0}^{n}(b_{i}a_{n-i}+c_{i}a_{n-i}) \\
		&=&\sum_{i=0}^{n}b_{i}a_{n-i}+\sum_{i=0}^{n}c_{i}a_{n-i}=(b\ast
		a)_{n}+(c\ast a)_{n}=((b\ast a)+(c\ast a))_{n}
	\end{eqnarray*}
	
	$\because $%
	\begin{equation*}
	((b+c)\ast a)_{n}=((b\ast a)+(c\ast a))_{n},\left( n\leq k\right) 
	\end{equation*}

	$\therefore $%
	\begin{equation*}
	(b+c)\ast a=a\ast (b+c),\ (a,b,c\in \mathcal{A}_{k}(R)
	\end{equation*}
	
	$\therefore $ \qquad $\mathbf{Right\_Distributivity}\left( \mathcal{A}_{k}(R),\ast
	,.\right) $
	
	\vspace{0.7em}
	\item[\protect\underline{Multiplicative Identity}] 
	\text{  \  }\\
	$\because $%
	\begin{equation*}
	(1\ast a)_{n}=\sum_{i=0}^{n}1_{i}a_{n-i}=1_{R}\cdot a_{n}=a_{n}
	\end{equation*}
	
	and%
	\begin{equation*}
	(a\ast 1)_{n}=\sum_{i=0}^{n}a_{i}1_{n-i}=a_{n}\cdot 1_{R}=a_{n}
	\end{equation*}
	
	$\therefore $ 
	\begin{equation*}
	1\ast a=a\ast 1=a
	\end{equation*}
	
	$\therefore $ $\mathbf{Multiplicative\_Identity}\left( \mathcal{A}%
	_{k}(R),\ast ,.\right) $
	\vspace{0.7em}
	\item[\protect\underline{Algebra Compatibility}] 
	\text{  \  }\\ %
	
	$\because $%
	\begin{eqnarray*}
		(r\cdot (a\ast b))_{n} &=&r(a\ast
		b)_{n}=r\sum_{i=0}^{n}a_{i}b_{n-i}=\sum_{i=0}^{n}(ra_{i})b_{n-i}=((r\cdot
		a)\ast b)_{n} \\
		&=&\sum_{i=0}^{n}a_{i}(rb_{n-i})=(a\ast (r\cdot b))_{n}
	\end{eqnarray*}
	
	$\therefore $ 
	\begin{equation*}
	r\cdot (a\ast b)=(r\cdot a)\ast b=a\ast (r\cdot b)
	\end{equation*}
	
	$\therefore \mathbf{Algebra\_Compatibility}\left( \mathcal{A}_{k}(R),\ast
	\right) $
	\vspace{0.7em}
	\item[\protect\underline{Commutativity when $R$ Commutative}] 
	\text{\ }
	\\
	
	$\because $%
	\begin{equation*}
	(a\ast
	b)_{n}=\sum_{i=0}^{n}a_{i}b_{n-i}=\sum_{j=0}^{n}a_{n-j}b_{j}=%
	\sum_{j=0}^{n}b_{j}a_{n-j}=(b\ast a)_{n}
	\end{equation*}
	
	$\therefore $ 
	\begin{equation*}
	a\ast b=b\ast a,\text{ \ \ \ }a,b\in \mathcal{A}_{k}(R)
	\end{equation*}
	
	$\therefore $%
	\begin{equation*}
	\mathbf{Commutative}\left( R,. \right) \Rightarrow \mathbf{Commutative}%
	\left( \mathcal{A}_{k}(R),\ast \right) 
	\end{equation*}
	
	\vspace{0.5in}
	
	\item[Summary] $(\mathcal{A}_{k}(R),+,\ast ,\cdot )$ associative $R$-algebra.
	
	Dimension $k+1$ as $R$-module.
	
	Identity: $1=(1_{R},0,\dots ,0)$.
	
	If $R$ commutative: $\mathcal{A}_{k}(R)$ commutative.
\end{description}
	\begin{theorem}
	If	\(R\) is an abstract commutative ring, \(u,v,t \in R\), \(n\in\mathbb{N}\), then the Binomial Chu-Vandermonde Identity \eqref{Ring-CV} in  this ring holds
	\begin{equation}
	(u+v)^{\,n:t}=\sum_{k=0}^{n}\binom{n}{k}u^{(n-k):t}v^{k:t},\label{Ring-CV}
	\end{equation}
	such that:
	\[
	u^{n:t} :=\prod_{j=0}^{n-1}(u+j t)\qquad(u^{0:t}:=1).
	\]
\end{theorem}

\begin{proof}
	\scriptsize
	\setlength{\LTleft}{-0.5in} 
	\setlength{\LTright}{0pt} 
	\begin{prooftable}
		1 & \((u+v)^{0:t}=1\) & (Def) \\
		\hline
		2 & \(\displaystyle\sum_{k=0}^{0}\binom{0}{k}u^{0:t}v^{0:t}=1\) & (Def, \(\binom{0}{0}=1\)) \\
		\hline
		3 & \(\displaystyle (u+v)^{n:t}=\sum_{k=0}^{n}\binom{n}{k}u^{(n-k):t}v^{k:t}\) & (IH) \\
		\hline
		4 & \((u+v)^{(n+1):t} = (u+v)^{n:t}\cdot(u+v+nt)\) & (Def) \\
		\hline
		5 & \((u+v)^{(n+1):t} = \left(\sum_{k=0}^{n}\binom{n}{k}u^{(n-k):t}v^{k:t}\right)(u+v+nt)\) & (4, 3) \\
		\hline
		6 & \(u+v+nt = (u+(n-k)t) + (v+kt)\) & (Alg) \\
		\hline
		7 & \((u+v)^{(n+1):t} = \sum_{k=0}^{n}\binom{n}{k}u^{(n-k):t}(u+(n-k)t)v^{k:t} + \sum_{k=0}^{n}\binom{n}{k}u^{(n-k):t}v^{k:t}(v+kt)\) & (Dist, 6) \\
		\hline
		8 & \(u^{(n-k):t}(u+(n-k)t) = u^{(n-k+1):t}\) & (Rec) \\
		\hline
		9 & \(v^{k:t}(v+kt) = v^{(k+1):t}\) & (Rec) \\
		\hline
		10 & \((u+v)^{(n+1):t} = \sum_{k=0}^{n}\binom{n}{k}u^{(n-k+1):t}v^{k:t} + \sum_{k=0}^{n}\binom{n}{k}u^{(n-k):t}v^{(k+1):t}\) & (8, 9, 7) \\
		\hline
		11 & \(\displaystyle\sum_{k=0}^{n}\binom{n}{k}u^{(n-k):t}v^{(k+1):t} = \sum_{j=1}^{n+1}\binom{n}{j-1}u^{(n-j+1):t}v^{j:t}\) & (Reindex: \(j=k+1\)) \\
		\hline
		12 & \((u+v)^{(n+1):t} = \sum_{k=0}^{n}\binom{n}{k}u^{(n-k+1):t}v^{k:t} + \sum_{k=1}^{n+1}\binom{n}{k-1}u^{(n-k+1):t}v^{k:t}\) & (10, 11) \\
		\hline
		13 & \(\displaystyle (u+v)^{(n+1):t} = \sum_{k=0}^{n+1}\binom{n+1}{k}u^{(n+1-k):t}v^{k:t}\) & (BinRec: \(\binom{n}{k}+\binom{n}{k-1}=\binom{n+1}{k}\)) \\
		\hline
		14 & \(\displaystyle (u+v)^{n:t}=\sum_{k=0}^{n}\binom{n}{k}u^{(n-k):t}v^{k:t}\) & (Induction) \\
	\end{prooftable}
\end{proof}
\section{Geometry}
\begin{definition}
	Set $p$ \ a prime, and $n\in 
	\mathbb{N}
	,$ $n\leq p.$
	
	Set $u,v$ polynumbers in basis $\left\{ e_{i}\right\} _{i=0,p-1}$ and $%
	\left\{ \alpha _{i}\right\} _{=0,p-1}$ $\subset $ $\mathbb{F}_{p}.$
	
	We define:
	
	\begin{itemize}
		\item $\mathcal{P}_{n}\left( u\right) \mathcal{ =}\left\{ \alpha
		.u^{k},\text{ }k\in 
		\mathbb{N}
		,\text{ }k\geq n,\text{ }\alpha \in \mathbb{F}_{p}\right\} $
		
		\item  $v\equiv 0\ \ \left( \ \mathrm{mod}\text{ }u^{n}\right)
		\Leftrightarrow \exists w\in \mathcal{P}_{n}\left( u\right) \mathcal{%
			 }$: $v=w.u$
		
		\item $u=\sum\limits_{i=0}^{n}\alpha _{i}e_{i},$ $\alpha _{n}\neq
		0\Rightarrow \mathbf{order}\left( u\right) =n$
	\end{itemize}
\end{definition}

\begin{corollary}

	{Set $p$}${>2}$ a prime{$,$ }$n\in 
	\mathbb{N}
	,$ $n\leq p,$ and $u$ a polynumber of order $n,$ then the equivalence %
	\eqref{QuatEquaEquipolyn} holds for any polynumber $x$ 
	\begin{equation}
	\left( x^{2}-x+u\equiv 0\ \ \left( \mathrm{mod}\text{ }u^{n}\right) \right)
	\Leftrightarrow \left( x\equiv \sum\limits_{k=1}^{n-1}C_{k-1}\cdot \,u^{k}\
	\ \left( \mathrm{mod}\text{ }u^{n}\right) \right)  \label{QuatEquaEquipolyn}
	\end{equation}
\end{corollary}
\begin{example}
	\text{\ \ }
	
	\begin{itemize}
		\item $p=3,$ $u\equiv e_{0}+e_{1},$ $\mathbf{order}\left( u\right) =1$
		
		$u_{1}=u^{1}\equiv e_{0}+e_{1}$
		
		$u_{2}=u^{2}\equiv u^{2}\equiv e_{0}-e_{1}$
		
		$u_{3}=u^{3}\equiv e_{0}$
		
		$u_{4}=u^{4}\equiv u$
		
		$x\equiv C_{0}u+C_{1}u^{2}+C_{2}u^{3}\equiv \left( e_{0}+e_{1}\right)
		+\left( e_{0}-e_{1}\right) +2\left( e_{0}\right) \equiv \allowbreak e_{0}$
		
		$x^{2}\equiv e_{0}$
		
		$x^{2}-x+u\equiv u\equiv u^{1}$%
		\begin{equation*}
		x^{2}-x+u\equiv 0\ \ \left( \mathrm{mod}\text{ }u^{1}\right)
		\end{equation*}%
		\item $p=5,$ $u\equiv e_{0}+e_{4},$ $\mathbf{order}\left( u\right) =4$
		
		$u_{1}\equiv e_{0}+e_{4}$
		
		$u_{2}=u^{2}\equiv \allowbreak e_{0}+2e_{4}$
		
		$u_{3}=u^{3}\equiv e_{0}+3e_{4}$
		
		$u_{4}=u^{4}\equiv e_{0}+4e_{4}$
		
		$u_{5}=u^{5}\equiv e_{0}$%
		\begin{equation*}
		\left\{ 
		\begin{array}{c}
		u_{1}\equiv e_{0}+e_{4} \\ 
		u_{2}\equiv e_{0}+2e_{4}%
		\end{array}%
		\right. \Leftrightarrow \left\{ 
		\begin{array}{c}
		u_{2}-u_{1}\equiv e_{4} \\ 
		2u_{1}-u_{2}\equiv e_{0}%
		\end{array}%
		\right.
		\end{equation*}
		
		$x\equiv C_{0}u+C_{1}u^{2}+C_{2}u^{3}+C_{3}u^{4}=\allowbreak 2e_{0}+e_{4}$
		
		$x^{2}\equiv \left( 2e_{0}+e_{4}\right) ^{2}\equiv \allowbreak -e_{0}-e_{4}$
		
		$x^{2}-x+u\equiv \allowbreak -2e_{0}-e_{4}\equiv -2\left( u_{2}-u_{1}\right)
		-\left( 2u_{1}-u_{2}\right) \equiv \allowbreak -u_{2}\equiv 4u^{7}\equiv
		4u^{3}.u^{4}$%
		\begin{equation*}
		x^{2}-x+u\equiv 0\ \ \left( \mathrm{mod}\text{ }u^{4}\right)
		\end{equation*}
	\end{itemize}
\end{example}

\begin{definition}
	For a field $K$ with characteristic $\mathrm{char}(K)\neq 2,3$ and
	parameters $a,b\in K$ with $4a^{3}+27b^{2}\neq 0$, define the short
	Weierstrass form of an elliptic curve over $K$ by:%
	\begin{equation*}
	E_{a,b}\left( K\right) :\qquad y^{2}=x^{3}+ax+b.
	\end{equation*}%
	The projective closure is%
	\begin{equation*}
	zy^{2}=x^{3}+axz^{2}+bz^{3}\subset \mathbb{P}_{K}^{2},
	\end{equation*}%
	with identity point (neutral element) $\mathcal{O}=[0:1:0]$.\newline
	For points $P=(x_{1},y_{1}),\;Q=(x_{2},y_{2})\in E_{a,b}\left( K\right) (K)$
	the following formulas hold.%
	\begin{equation*}
	-P=(x_{1},-y_{1}).
	\end{equation*}%
	\begin{equation*}
	P+Q=\left\{ 
	\begin{array}{lcc}
	\left( \lambda ^{2}-x_{1}-x_{2},\lambda (x_{1}-x_{3})-y_{1}\right) &  & 
	x_{1}\neq x_{2}\text{, \ }\lambda =\frac{y_{2}-y_{1}}{x_{2}-x_{1}} \\ 
	\mathcal{O} &  & x_{1}=x_{2},\, y_{1}\neq y_{2}%
	\end{array}%
	\right.
	\end{equation*}%
	\begin{equation*}
	2P=\left\{ 
	\begin{array}{lcc}
	\left( \lambda ^{2}-2x_{1},\lambda (x_{1}-x_{3})-y_{1}\right) &  & y_{1}\neq
	0,\text{ \ }\lambda =\frac{3x_{1}^{2}+a}{2y_{1}} \\ 
	\mathcal{O} &  & y_{1}=0%
	\end{array}%
	\right.
	\end{equation*}
\end{definition}
\text{\ }

\begin{example}
	Let $K=\mathbb{F}_{5}$ and consider the elliptic curve in short Weierstrass
	form%
	\begin{equation*}
	E_{1,1}:\qquad y^{2}=f(x)=x^{3}+x+1.
	\end{equation*}%
	For a curve $y^{2}=f_{a,b}(x)=x^{3}+ax+b$ the discriminant condition nonzero
	is detected by%
	\begin{equation*}
	\Delta \coloneqq-16(4a^{3}+27b^{2}),
	\end{equation*}%
	so it suffices to check $4a^{3}+27b^{2}\not\equiv 0\pmod{5}$. With $a=b=1$
	we have%
	\begin{equation*}
	4\cdot 1^{3}+27\cdot 1^{2}=4+27=31\equiv 1\pmod{5},
	\end{equation*}%
	hence $\Delta \neq 0$ and $E_{1,1}$ is nonsingular over $\mathbb{F}_{5}$.
	Compute the squares:%
	\begin{equation*}
	0^{2}\equiv 0,\quad 1^{2}\equiv 1,\quad 2^{2}\equiv 4,\quad 3^{2}\equiv
	9\equiv 4,\quad 4^{2}\equiv 16\equiv 1\pmod{5}.
	\end{equation*}%
	Thus the set of quadratic residues (squares) in $\mathbb{F}_{5}$ is $%
	\{0,1,4\}. $ 
	
	\text{\  }
	
	\bigskip
	
	\text{\  }
	
	\begin{align*}
	f(0)& =0^{3}+0+1=1\in \{0,1,4\}\Rightarrow y^{2}=1\Rightarrow y\in \left\{
	1,-1\right\} _{5}\equiv \left\{ 1,4\right\} _{5} \\
	f(1)& =1^{3}+1+1=3\notin \{0,1,4\}\Rightarrow y^{2}=3\text{ has no solution
		in }\mathbb{F}_{5}, \\
	f(2)& =2^{3}+2+1=8+2+1=11\equiv 1\pmod{5}\Rightarrow y\in \left\{
	1,4\right\} _{5}, \\
	f(3)& =3^{3}+3+1=27+3+1=31\equiv 1\pmod{5}\Rightarrow y\in \left\{
	1,4\right\} _{5}, \\
	f(4)& =4^{3}+4+1\equiv 4\pmod{5}\Rightarrow y^{2}=4\Rightarrow y\in \left\{
	2,-2\right\} _{5}\equiv \left\{ 2,3\right\} _{5}
	\end{align*}%
	\begin{equation*}
	E_{1,1}(\mathbb{F}_{5})=\left\{ \,\mathcal{O}%
	,\,(0,1),(0,4),(2,1),(2,4),(3,1),(3,4),(4,2),(4,3)\,\right\} .
	\end{equation*}%
	\begin{equation*}
	\#E_{1,1}(\mathbb{F}_{5})=8+1=9.
	\end{equation*}%
	\begin{table}[H]
		\begin{tabular}{|c|c|c|c|c|}
			\hline
			\textbf{k} & \textbf{Operation} & $\boldsymbol{\lambda }$ & $\boldsymbol{%
				x_{k}}$ & $\boldsymbol{y_{k}}$ \\ \hline
			1 & \multicolumn{1}{|l|}{$P=(0,1)$} & -- & $0$ & $1$ \\ \hline
			2 & \multicolumn{1}{|l|}{$2P=P+P$} & $\frac{3x_{1}^{2}+a}{2y_{1}}=\frac{1}{%
				2}\equiv 3$ & $4$ & $2$ \\ \hline
			3 & \multicolumn{1}{|l|}{3P= $2P+P=(4,2)+(0,1)$} & $\frac{1-2}{0-4}\equiv 4$
			& $2$ & $1$ \\ \hline
			4 & \multicolumn{1}{|l|}{$4P=2P+2P=(4,2)$} & $\frac{3\cdot 4^{2}+a}{2\cdot 2%
			}\equiv 1$ & $3$ & $4$ \\ \hline
			5 & \multicolumn{1}{|l|}{$5P= 4P+P=(3,4)+(0,1)$} & $\frac{1-4}{0-3}\equiv 1$
			& $3$ & $1$ \\ \hline
			6 & \multicolumn{1}{|l|}{$6P=3P+3P=(2,1)$} & $\frac{3\cdot 2^{2}+a}{2\cdot 1%
			}\equiv 4$ & $2$ & $4$ \\ \hline
			7 & \multicolumn{1}{|l|}{$7P= 6P+P=(2,4)+(0,1)$} & $\frac{1-4}{0-2}\equiv 4$
			& $4$ & $3$ \\ \hline
			8 & \multicolumn{1}{|l|}{$8P= 7P+P=(4,3)+(0,1)$} & $\frac{1-3}{0-4}\equiv 3$
			& $0$ & $4$ \\ \hline
			9 & \multicolumn{1}{|l|}{$9P= 8P+P=(0,4)+(0,1)$} & -- & $\mathcal{O}$ & --
			\\ \hline
		\end{tabular}%
		
	\end{table}
	\begin{table}[H]
		\begin{equation*}
		\begin{aligned} &0P =\mathcal{O},\qquad 1P = (0,1),\qquad 2P = (4,2),\\ &3P
		= (2,1),\qquad 4P = (3,4),\qquad 5P = (3,1),\\ &6P = (2,4),\qquad 7P =
		(4,3),\qquad 8P = (0,4). \end{aligned}
		\end{equation*}
		\caption{ Cyclic ordering of points on $E_{1,1}$ over $\mathbb{F}_{5}$
			(generator $P=(0,1)$).}
	\end{table}
	\begin{table}[H]
		\begin{equation*}
		\renewcommand{\arraystretch}{1.05} 
		\begin{array}{c|ccccccccc}
		+ & 0P & 1P & 2P & 3P & 4P & 5P & 6P & 7P & 8P \\ \hline
		0P & \mathcal{O} & (0,1) & (4,2) & (2,1) & (3,4) & (3,1) & (2,4) & (4,3) & 
		(0,4) \\ 
		1P & (0,1) & (4,2) & (2,1) & (3,4) & (3,1) & (2,4) & (4,3) & (0,4) & 
		\mathcal{O} \\ 
		2P & (4,2) & (2,1) & (3,4) & (3,1) & (2,4) & (4,3) & (0,4) & \mathcal{O} & 
		(0,1) \\ 
		3P & (2,1) & (3,4) & (3,1) & (2,4) & (4,3) & (0,4) & \mathcal{O} & (0,1) & 
		(4,2) \\ 
		4P & (3,4) & (3,1) & (2,4) & (4,3) & (0,4) & \mathcal{O} & (0,1) & (4,2) & 
		(2,1) \\ 
		5P & (3,1) & (2,4) & (4,3) & (0,4) & \mathcal{O} & (0,1) & (4,2) & (2,1) & 
		(3,4) \\ 
		6P & (2,4) & (4,3) & (0,4) & \mathcal{O} & (0,1) & (4,2) & (2,1) & (3,4) & 
		(3,1) \\ 
		7P & (4,3) & (0,4) & \mathcal{O} & (0,1) & (4,2) & (2,1) & (3,4) & (3,1) & 
		(2,4) \\ 
		8P & (0,4) & \mathcal{O} & (0,1) & (4,2) & (2,1) & (3,4) & (3,1) & (2,4) & 
		(4,3)%
		\end{array}%
		\end{equation*}%
		\caption{ Cayley addition table for $E_{1,1}(\mathbb{F}_{5})$ in point
			notation.}
	\end{table}
	\begin{table}[H]
		\begin{equation*}
		\begin{array}{c|ccccccccc}
		+ & 0 & 1 & 2 & 3 & 4 & 5 & 6 & 7 & 8 \\ \hline
		0 & 0 & 1 & 2 & 3 & 4 & 5 & 6 & 7 & 8 \\ 
		1 & 1 & 2 & 3 & 4 & 5 & 6 & 7 & 8 & 0 \\ 
		2 & 2 & 3 & 4 & 5 & 6 & 7 & 8 & 0 & 1 \\ 
		3 & 3 & 4 & 5 & 6 & 7 & 8 & 0 & 1 & 2 \\ 
		4 & 4 & 5 & 6 & 7 & 8 & 0 & 1 & 2 & 3 \\ 
		5 & 5 & 6 & 7 & 8 & 0 & 1 & 2 & 3 & 4 \\ 
		6 & 6 & 7 & 8 & 0 & 1 & 2 & 3 & 4 & 5 \\ 
		7 & 7 & 8 & 0 & 1 & 2 & 3 & 4 & 5 & 6 \\ 
		8 & 8 & 0 & 1 & 2 & 3 & 4 & 5 & 6 & 7%
		\end{array}%
		\end{equation*}%
		\caption{ Addition modulo $9$ corresponding to  the cyclic group generated by $%
			P$.}
	\end{table}
	
	\begin{figure}[H]
		\centering
		\begin{tikzpicture}
		\begin{axis}[
		title={Plot of $y^{2}=x^{3}+x+1$ over $\mathbb{F}_{5}$ },
		width=10cm, height=8cm,
		xmin=-0.5, xmax=4.5, ymin=-0.5, ymax=4.5,
		xtick={0,1,2,3,4}, ytick={0,1,2,3,4},
		axis lines=middle, xlabel=$x$, ylabel=$y$, grid=both,
		minor tick num=0
		]
		\addplot [very thick, black, domain=-0.5:4.5, samples=800, smooth] 
		({x}, {sqrt(abs(x^3 + x + 1))}) ;
		\addplot [very thick, black, domain=-0.5:4.5, samples=800, smooth] 
		({x}, {-sqrt(abs(x^3 + x + 1))}) ;
		
		\addplot[only marks, mark=*, mark size=3pt, black] coordinates {
			(0,1) (0,4) (2,1) (2,4) (3,1) (3,4) (4,2) (4,3)
		};
		
		\node at (axis cs:0.15,1.15) {\small $(0,1)$};
		\node at (axis cs:0.15,3.85) {\small $(0,4)$};
		\node at (axis cs:2.15,1.15) {\small $(2,1)$};
		\node at (axis cs:2.15,3.85) {\small $(2,4)$};
		\node at (axis cs:3.15,1.15) {\small $(3,1)$};
		\node at (axis cs:3.15,3.85) {\small $(3,4)$};
		\node at (axis cs:4.15,2.15) {\small $(4,2)$};
		\node at (axis cs:4.15,2.85) {\small $(4,3)$};
		
		\node[draw,inner sep=1pt,circle] at (axis cs:4.0,4.2) {$O$};
		\end{axis}
		\end{tikzpicture}
	\end{figure}
\end{example}

\makeatletter\let\orig@addcontentsline%
\addcontentsline
\renewcommand{\addcontentsline}[3]{}   \makeatother

\makeatletter\let\addcontentsline\orig@addcontentsline \makeatother
\end{document}